\documentclass[letterpaper,reqno,11pt,oneside,final]{amsart}
\pdfoutput=1
\usepackage{amsmath,amsthm,amsfonts,amssymb}
\usepackage[extension=pdf,bookmarksopen]{hyperref}
\usepackage{enumitem,comment}
\usepackage[T1]{fontenc}
\usepackage{times}
\usepackage{mathtools}
\usepackage{mathrsfs}
\usepackage{comment}
\usepackage{cases}
\usepackage{bm}
\usepackage[dvipsnames]{xcolor}
\usepackage[multiple]{footmisc}
\usepackage{microtype}
\usepackage[totalheight=9.6in,totalwidth=6.15in,centering]{geometry}
\usepackage{graphics,graphicx}
\usepackage{tocvsec2}
\usepackage[color,notref,notcite]{showkeys}
\usepackage[style=alphabetic,sorting=nyt,natbib=true,maxnames=99,isbn=false,doi=false,url=false,firstinits=true,hyperref=auto,arxiv=abs,backend=bibtex]{biblatex}
\AtEveryBibitem{%
  \clearlist{language}}
\DeclareFieldFormat[article,inbook,incollection,inproceedings,patent,thesis,unpublished]{title}{#1\isdot}
\renewbibmacro{in:}{%
  \ifentrytype{article}{}{%
    \printtext{\bibstring{in}\intitlepunct}}}
\defbibheading{apa}[\refname]{\section*{#1}}
\DeclareFieldFormat{sentencecase}{\MakeSentenceCase{#1}} \renewbibmacro*{title}{%
  \ifthenelse{\iffieldundef{title}\AND\iffieldundef{subtitle}}{}
  {\ifthenelse{\ifentrytype{article}\OR\ifentrytype{inbook}%
      \OR\ifentrytype{incollection}\OR\ifentrytype{inproceedings}%
      \OR\ifentrytype{inreference}} {\printtext[title]{%
        \printfield[sentencecase]{title}%
        \setunit{\subtitlepunct}%
        \printfield[sentencecase]{subtitle}}}%
    {\printtext[title]{%
        \printfield[titlecase]{title}%
        \setunit{\subtitlepunct}%
        \printfield[titlecase]{subtitle}}}%
    \newunit}%
  \printfield{titleaddon}}

\definecolor{darkblue}{rgb}{0.13,0.13,0.39}
\hypersetup{colorlinks=true,urlcolor=darkblue,citecolor=darkblue,linkcolor=darkblue,%
  pdftitle=KP(Z),pdfauthor={J. Quastel, D. Remenik}}

\DeclareSymbolFont{eulerletters}{U}{eur}{m}{n}

\mathtoolsset{showmanualtags,showonlyrefs}

\newtheorem{thm}{Theorem}[section]

\theoremstyle{definition}
\newtheorem{rem}[thm]{Remark}

\newtheorem{ex}[thm]{Example}

\newcommand{\fh}{\mathfrak{h}}

\newcommand{\R}{R}

\newcommand{\I}{\uptext{i}}
\newcommand{\pp}{\mathbb{P}}

\newcommand{\ee}{\mathbb{E}}
\newcommand{\rr}{\mathbb{R}}

\newcommand{\aip}{\mathcal{A}}
\newcommand{\aipo}{\mathcal{A}_1}

\newcommand{\cD}{\mathcal{D}}

\newcommand{\p}{\partial}
\newcommand{\uno}[1]{\mathbf{1}_{#1}}
\newcommand{\ep}{\varepsilon}

\newcommand{\vs}{\vspace{6pt}}
\newcommand{\wt}{\widetilde}
\newcommand{\wh}{\widehat}
\newcommand{\K}{K_{\Ai}}
\newcommand{\m}{\overline{m}}

\newcommand{\qand}{\quad\text{and}\quad}
\newcommand{\qqand}{\qquad\text{and}\qquad}

\newcommand{\fU}{\mathbf{U}}

\newcommand{\D}{\mathrm{d}}
\renewcommand{\d}{\mathrm{d}}
\newcommand{\fR}{\mathbf{R}}
\newcommand{\q}{\mathfrak{q}}

\makeatletter
\newcommand{\na}[1][\@nil]{%
  \def\tmp{#1}%
   \ifx\tmp\@nnil
       \nabla_{\!a}\tts
    \else
         \nabla^{#1}_{\!a}\tts
    \fi}
\makeatother

\newcommand{\ts}{\hspace{0.1em}}
\newcommand{\tts}{\hspace{0.05em}}
\newcommand{\tsm}{\hspace{-0.1em}}
\newcommand{\ttsm}{\hspace{-0.05em}}

\newcommand{\itwopii}[1]{\frac{1}{(2\pi\I)^{#1}}}

\DeclareMathOperator\arctanh{arctanh}

\makeatletter
\newcommand\RedeclareMathOperator{%
  \@ifstar{\def\rmo@s{m}\rmo@redeclare}{\def\rmo@s{o}\rmo@redeclare}%
}
\newcommand\rmo@redeclare[2]{%
  \begingroup \escapechar\m@ne\xdef\@gtempa{{\string#1}}\endgroup
  \expandafter\@ifundefined\@gtempa
     {\@latex@error{\noexpand#1undefined}\@ehc}%
     \relax
  \expandafter\rmo@declmathop\rmo@s{#1}{#2}}
\newcommand\rmo@declmathop[3]{%
  \DeclareRobustCommand{#2}{\qopname\newmcodes@#1{#3}}%
}
\@onlypreamble\RedeclareMathOperator
\makeatother
\newcommand{\uptext}[1]{\text{\upshape{#1}}}

\DeclareMathOperator{\hypo}{\uptext{hypo}}

\DeclareMathOperator{\UC}{\uptext{UC}}

\DeclareMathOperator{\Ai}{\uptext{Ai}}

\DeclareMathOperator{\tr}{\uptext{tr}}
\RedeclareMathOperator{\det}{\mathop{\uptext{det}}}
\RedeclareMathOperator{\ker}{\mathop{\uptext{ker}}}
\RedeclareMathOperator{\exp}{\mathop{\uptext{exp}}}
\RedeclareMathOperator{\log}{\mathop{\uptext{log}}}
\RedeclareMathOperator*{\lim}{\mathop{\uptext{lim}}}
\RedeclareMathOperator*{\sup}{\mathop{\uptext{sup}}}
\RedeclareMathOperator*{\limsup}{\mathop{\uptext{lim\hspace{1pt}sup}}}
\RedeclareMathOperator*{\liminf}{\mathop{\uptext{lim\hspace{1pt}inf}}}
\RedeclareMathOperator*{\max}{\mathop{\uptext{max}}}
\RedeclareMathOperator*{\inf}{\mathop{\uptext{inf}}}
\RedeclareMathOperator*{\arctanh}{\mathop{\uptext{arctanh}}}
\RedeclareMathOperator*{\min}{\mathop{\uptext{min}}}
\RedeclareMathOperator*{\cos}{\mathop{\uptext{cos}}}
\RedeclareMathOperator*{\sin}{\mathop{\uptext{sin}}}
\RedeclareMathOperator*{\arg}{\mathop{\uptext{arg}}}
\RedeclareMathOperator{\Re}{\uptext{Re}}
\RedeclareMathOperator{\Im}{\uptext{Im}}

\newcommand{\fT}{\mathbf{S}}
\newcommand{\ft}{\mathbf{t}}
\newcommand{\fs}{\mathbf{s}}
\newcommand{\s}{\fs}
\renewcommand{\b}{\mathbf{b}}
\newcommand{\fx}{\mathbf{x}}

\newcommand{\fA}{\mathbf{A}}

\newcommand{\fB}{B}
\newcommand{\fBb}{\mathbf{B}}
\newcommand{\fP}{\mathbf{P}}
\newcommand{\fL}{\mathbf{L}}
\newcommand{\fK}{\mathbf{K}}

\newcommand{\fM}{\mathbf{M}}

\newcommand{\fI}{\mathbf{I}}

\newcommand{\g}{\bar{\bm{\gamma}}}

\renewcommand{\P}{\chi}
\newcommand{\kh}{{\rm h}}

\def\dash---{\kern.16667em---\penalty\exhyphenpenalty\hskip.16667em\relax}

\numberwithin{equation}{section}

\let\oldmarginpar\marginpar
\renewcommand\marginpar[1]{\-\oldmarginpar[\raggedleft\footnotesize #1]%
  {\raggedright{\small\textsf{#1}}}}

\allowdisplaybreaks[2]

\begin{document}

\maxtocdepth{subsection}

\title[KP(Z)]{KP governs random growth off a one dimensional substrate}

\date{August 27, 2021}

\author{Jeremy Quastel} \address[J.~Quastel]{
  Department of Mathematics\\
  University of Toronto\\
  40 St. George Street\\
  Toronto, Ontario\\
  Canada M5S 2E4} \email{quastel@math.toronto.edu}

\author{Daniel Remenik} \address[D.~Remenik]{
  Departamento de Ingenier\'ia Matem\'atica and Centro de Modelamiento Matem\'atico (UMI-CNRS 2807)\\
  Universidad de Chile\\
  Av. Beauchef 851, Torre Norte, Piso 5\\
  Santiago\\
  Chile} \email{dremenik@dim.uchile.cl}

\begin{abstract}
The logarithmic derivative of the marginal distributions of randomly fluctuating interfaces in one dimension on a large scale evolve according to the Kadomtsev--Petviashvili (KP) equation.
This is derived algebraically from a Fredholm determinant obtained in \cite{fixedpt} for the KPZ fixed point as the limit of the transition probabilities of TASEP, a special  solvable model in the KPZ universality class.
The Tracy-Widom distributions appear as special self-similar solutions of KP and KdV.
In addition, it is noted that several known exact solutions of the KPZ equation also solve KP.
\end{abstract}

\dedicatory{Dedicated to the memory of Harold Widom}

\maketitle
\tableofcontents

\section{Matrix KP equation for  multidimensional distributions}

The one dimensional  KPZ universality class consists of random growth models, last passage percolation and directed polymers, and random stirred fluids.  All models in the class have an analogue of the height function $h(t,x)$ (free energy, integrated velocity) whose
long time large scale evolution is the principal object of study.
The name of the class comes from the Kardar-Parisi-Zhang equation,
\begin{equation}\label{KPZ}
\partial_t h = \lambda(\partial_xh)^2  + \nu \partial_x^2h + \sigma \xi
\end{equation}
with $\xi$ a space-time white noise, a canonical continuum equation for random growth introduced in \cite{kpz}.
However, the key interest is on the universal features which are only found in large space-time scales, under the 1:2:3 scaling corresponding to $\ep\to0$ in
\begin{equation}\label{eq:123sc}
\ep^{1/2} h(\ep^{-3/2} t, \ep^{-1} x) - C_\ep t. 
\end{equation}
The KPZ equation is not invariant under this scaling, which sends $(\lambda,\nu,\sigma)$ to $(\lambda,\ep^{1/2}\nu,\ep^{1/4}\sigma)$.  A key problem is to find the true, scaling invariant equation for random interface growth.

Since the early 2000's \cite{Johansson,sasamoto,borFerPrahSasam} it was known, for a number of models in the class, and special scaling invariant initial data \emph{narrow wedge} and \emph{flat}, that the distributional limits of \eqref{eq:123sc} were the Tracy-Widom distributions of random matrix theory. In an earlier article \cite{fixedpt} it was shown that, at least for one model in the class, TASEP, \eqref{eq:123sc} converges to a Markov process $\fh(t,x)$ which is invariant under 1:2:3 scaling,
\begin{equation}\label{eq:123inv}
\alpha\fh(\alpha^{-3}t,\alpha^{-2}x;\alpha^{-1}\fh_0(\alpha^2x))\stackrel{\uptext{dist}}{=}\fh(t,x;\fh_0),
\end{equation}
where $\fh_0$ after the semicolon denotes the initial data.
It is widely believed that this \emph{KPZ fixed point} governs the limiting fluctuation for all models in the class.

The KPZ fixed point does not satisfy a stochastic differential equation.
In place of that, it inherits a variational formulation from TASEP; a Hopf-Lax type formula involving a non-trivial input noise called the \emph{Airy sheet} $\aip(x,y)$: for the KPZ fixed point starting from $\fh(0,x)=\fh_0(x)$,
\begin{equation}\label{eq:var}
\fh(t,x)  \stackrel{\uptext{dist}}{=}  \sup_{y\in\rr}\big\{ t^{1/3}\aip(t^{-2/3} x,t^{-2/3} y)- \tfrac1{t}(x-y)^2 + \fh_0(y)\big\}.
\end{equation}
The Airy sheet $\aip(x,y)$ can be thought of as the height function at $x$ at time $1$, starting from a narrow wedge at $y$ at time $0$, and therefore involves coupling different initial conditions.
As far as we know at the present time, the coupled initial condition problem is not integrable, and therefore the distribution of the Airy sheet is unknown.
This led to a problem in that it was unclear that \eqref{eq:var} even involved a unique object on the right hand side.  An important advance is in \cite{DOV}, who show
that the Airy sheet is a functional of the Airy line ensemble. 
This puts the variational formula \eqref{eq:var} on a solid footing, as it obviates the need for uniqueness of the Airy sheet.
However, the functional is completely non-explicit. 
In this sense, \eqref{eq:var} is not satisfying as a universal scaling invariant equation.

Instead of a universal stochastic equation, one can study the $n$-space point distribution functions,
\begin{equation}\label{eq:F1}
F(t,x_1,\dotsc,x_n,r_1,\dotsc, r_n) = \pp_{\fh_0}\!\left(\fh(t,x_1)\leq r_1,\dotsc,\fh(t,x_n)\leq r_n\right)
\end{equation}
where the subscript in $\pp_{\fh_0}$ denotes the KPZ fixed point initial data.
In the cases of narrow wedge and flat initial data, it was known \cite{johanssonShape,sasamoto,borFerPrahSasam} that the one-dimensional distributions $F(1, x,r)$ were, respectively, the Tracy-Widom GUE and GOE random matrix distributions (but, except in the particular case of narrow wedge initial data, the connection between random growth and random matrices has remained tangential and murky). 
The multidimensional distributions in these cases are given by Fredholm determinants, and define the Airy$_2$ and Airy$_1$ processes.
The one-dimensional distributions can be written in terms of the Hastings-McLeod solution of Painlev\'e II; a longstanding open question was whether the distributions satisfy an equation in the more general setting.

\subsection{Main results}\label{sec:main}

In \cite{fixedpt}, it is shown that for initial data $\fh_0$ in $\UC$, a space of upper semicontinuous functions with a linear growth condition (see \eqref{eq:UC}), the $n$-space point distribution functions from \eqref{eq:F1} are given by Fredholm determinants 
\begin{equation}\label{eq:F1a}
F(t,x_1,\dotsc,x_n,r_1,\dotsc, r_n) = \det(\fI - \fK)
\end{equation}
where $\fK = \fK( t, x_1,\dotsc,x_n,r_1,\dotsc, r_n,\fh_0)$ is an operator on the $n$-fold direct sum of $L^2(\rr^+)$, given by an explicit $n\times n$ matrix kernel $\fK_{ij}(u_1,u_2)=\fK_{ij}(t, x_1,\dotsc,x_n,r_1,\dotsc, r_n,\fh_0,u_1,u_2)$, $i,j=1,\dotsc,n$.
The determinant is non-zero, so we can define $\fR= (\fI-\fK)^{-1}$.
Furthermore $\fK\fR= \fR\fK=\fI-\fR$ and we can let
\begin{equation}
\label{eq:F1b}
Q=\fR\fK(0,0),
\end{equation}
which is an $n\times n$ matrix valued function of $t,x_1,\dotsc,x_n,,r_1,\dotsc, r_n,$ and the initial height profile $\fh_0$.

Let
\begin{equation}\label{eq:drdx}
\cD_r= \partial_{r_1}+\cdots + \partial_{r_n}, \qquad \cD_x= \partial_{x_1}+\cdots + \partial_{x_n}.
\end{equation}
In Sec. \ref{sec:thirdsec} we will show that \begin{equation}\label{traceQ}
 \cD_r\log F = \tr Q.
\end{equation}
Our main result is: 

\begin{thm}\label{thm:2}
For any $\fh_0\in\UC$, $Q$ and its derivative $q= \cD_rQ$ solve the matrix Kadomtsev--Petviashvili (KP) equation
\begin{equation}\label{eq:matKP}
\p_tq+\tfrac12\tts\cD_rq^2+\tfrac1{12}\cD_r^3q+\tfrac14\cD_x^2Q+\tfrac12[q,\tts\cD_xQ]=0,
\end{equation}
where $[A,B] = AB-BA$.
In particular, for the one point marginals of the KPZ fixed point (i.e. \eqref{eq:F1} in  the case $n=1$),  $\phi=\partial_r^2 \log F$ satisfies the scalar KP-II equation
\begin{equation}\label{eq:KP-II}
\partial_t \phi + \tfrac12\partial_r \phi^2 + \tfrac1{12}\partial_r^3 \phi + \tfrac14\partial_r^{-1} \partial_x^2 \phi = 0.
\end{equation}
\end{thm}

The KP equation \eqref{eq:KP-II} was originally derived from studies of long waves in shallow water \cite{ablowitzSegur}.
It has come to be accepted as the natural two dimensional extension of the Korteweg--de Vries equation (KdV); when $\phi$ is independent of $x$, corresponding in our case to \emph{flat initial data}, it reduces to KdV,
\begin{equation}\label{eq:KdV}
\partial_t \phi + \tfrac12\partial_r \phi^2+ \tfrac1{12}\partial_r^3 \phi=0.
\end{equation}
KP is completely integrable and plays an important role in the
Sato theory as the first equation in the KP hierarchy \cite{jimboMiwaDate}.  The matrix KP equation \eqref{eq:matKP} exists in the literature, see e.g. \cite{koponel,sakhnovich}.
None of the previous physical derivations of KP seem to be related to the problem at hand, and it could well be that our evolution is through a 
class of functions where the equation is formally the same, because of similarities in the weakly nonlinear asymptotics, but the physics is completely different.  Note that the equations are usually written with with coefficients $3$, $1$ and $3$ replacing our $\tfrac12$, $\tfrac1{12}$ and $\tfrac14$, which is achieved by $t\mapsto 12t$ and $\phi\mapsto\tfrac12 \phi$.

\begin{rem}\label{rem:m}
\leavevmode
\begin{enumerate}[label=\arabic{*}.,itemsep=2pt]

\item The fact that random interface growth is governed by KP was not anticipated\footnote{However, see \cite{prolhacRiemann}, which appeared on the arXiv two days before this article was first posted, where it is shown that particular finite volume solutions \cite{baikLiu} can be written as superpositions of solitons. \cite{baikLiuSilva} treat other finite volume initial conditions.}.
We do not have physical intuition why it is true; it follows by, essentially, algebra from the form of the kernel in the Fredholm determinant for \eqref{eq:F1}, and we believe it is  the first example of a physical law having been obtained in such a fashion.
In retrospect, in the scalar case, it was known that the evolution equations (1)-(3) in Thm. \ref{thm3} below for the kernel lead to \eqref{eq:KP-II}, a fact that seems to have been rediscovered many times, particularly by \citet{poppeIP}, see also \cite{zaharovShabat2,poppePhysicaD,poppeSattinger,mckean-Fredholm}, before we rediscovered it.
The fact that the kernel evolves by Thm. \ref{thm3} (1)--(3), and its importance, was not recognized earlier and is one of the key contributions of this article.

\item There are not so many natural partial differential equations with the necessary (in view of \eqref{eq:123inv}) invariance under
\begin{equation}
\phi(t,x,r) \mapsto \alpha^{-2} \phi(\alpha^{-3} t,\alpha^{-2} x,\alpha^{-1} r),\qquad \fh_0(x)\mapsto \alpha^{-1} \fh_0(\alpha^2 x).
\end{equation}
So in retrospect, once one knows that the finite dimensional marginals come from a closed equation and that special examples connected to Painlev\'e transcendents, one might expect to look for integrable equations such as KP.\\
\emph{The question is why the finite dimensional marginals should come from a closed equation at all.} 

\item The \emph{Lax pair} formulation of \eqref{eq:KP-II} (and also 
\eqref{eq:matKP} with the matrix $q$ replacing $\phi$) is 
\begin{equation}\label{lax}
\partial_t L = [L,A], \qquad L= \partial_x + \partial^2_{r} + 2\phi.
\end{equation}
and $A= \tfrac13\partial_r^3 + \tfrac23 \phi \partial_r + \tfrac12 \partial_r\phi - \tfrac12 
\partial_r^{-1} \partial_x \phi$, which tells us that the ``spectrum'' of $L$ is conserved.  A very interesting question (a version of which A. Borodin asked us) is to make this precise and understand its physical meaning. 

\item The one dimensional distribution functions themselves  therefore satisfy the equivalent \emph{Hirota bilinear equation},
\begin{equation}
F\partial^2_{tr} F -\partial_t F\partial_r F+ \tfrac1{12}F\partial^4_r F - \tfrac13\partial_r F\partial^3_r F+ \tfrac14(\partial_r^2 F)^2+\tfrac14F\partial^2_x F-\tfrac14(\partial_x F)^2=0,
\end{equation}
which again has the necessary 1:2:3 invariance, now under
\begin{equation}
F(t,x,r) \mapsto F(\alpha^{-3} t,\alpha^{-2} x,\alpha^{-1} r),\qquad \fh_0(x)\mapsto \alpha^{-1} \fh_0(\alpha^2 x).
\end{equation}
\item Unlike other limit points for fluctuation universality classes in probability, the Tracy-Widom distributions themselves lack any invariance.  Thm. \ref{thm:2} recovers the invariance of the scaling limit under the 1:2:3 scaling, and the Tracy-Widom distributions then appear in the context of the KPZ universality class as special self-similar solutions of KP (see Sec. \ref{sec:TW}). 
\end{enumerate}
\end{rem}

Thm. \ref{thm:2} follows from the form of the kernel in the deteminantal formula \eqref{eq:F1a} and the following result.
In the scalar case this appeared earlier in \cite{poppeIP}.

\begin{thm}\label{thm3}
Let $\fK=\fK(t,x,r)$ be an operator on the $n$-fold direct sum of $L^2(\rr_{\ge 0})$ which is trace class uniformly in compact sets of $t>0$, $x,r\in\rr$, with matrix kernel $\fK(u,v)=\fK(t,x,r,u,v))$ which satisfies
\begin{enumerate}[label=\uptext{(\arabic*)}]
\item $\p_r\fK(u,v)=(\p_{u}+\p_{v})\fK$,
\item $\p_t\fK(u,v)=-\tfrac13(\p_{u}^3+\p_{v}^3)\fK(u,v)$,
\item $\p_x\fK(u,v)=(\p_{v}^2-\p_{u}^2)\fK(u,v)$.
\end{enumerate}
Suppose in addition that $\det(\fI-\fK) >0$ for all finite $t,x,r$, that $\fK$ is real analytic in $t$ and in each $x_i$ and $r_i$, and that the trace norm $\|\fK\|_1<1$ for $r$ in some open real interval.
Then the $n\times n$ matrix $Q=(\fI-\fK)^{-1}\fK(0,0)$  and its derivative  $q= \p_rQ$ satisfy the matrix Kadomtsev--Petviashvili (KP) equation
\begin{equation}\label{eq:matKPag}
\p_tq+\tfrac12\tts\p_rq^2+\tfrac1{12}\p_r^3q+\tfrac14\p_x^2Q+\tfrac12[q,\tts\p_xQ]=0.\end{equation}
\end{thm}

The kernel in \eqref{eq:F1a} is given explicitly (see \cite[Eqn. (4.2)]{fixedpt}) as
\begin{align} \label{eq:trcldecomp2}
&\fK_{ij}(t, x_1,\dotsc,x_n,r_1,\dotsc, r_n,\fh_0, u,v)= -e^{(x_j-x_i)\p^2}(u+r_i,v+r_j)\uno{x_i<x_j}\\&\qquad+
\int_{s\in\rr^+\atop z,b\in\rr } p_-(z,ds,db) \fT_{t,-x_i-s} (b,u+r_i)\fT_{t,x_j}(z,v+r_j) \\
&\qquad  + \int_{s\in\rr^+\atop z,b\in\rr } p_+(z,ds,db) \fT_{t,-x_i}  (z,u+r_i)  \fT_{t,x_j-s}  (b, v+r_j)\\
& \qquad - \int_{s_-,s_+\in\rr^+\atop z,b_-,b_+\in\rr } p_-(z,ds_-,db_-)p_+(z,ds_+,db_+) \fT_{t,-x_i-s_-}  (b_-,u+r_i)  \fT_{t,x_j-s_+}  (b_+, v+r_j),
\end{align} 
where $p_-(z,ds,db) = p_z(\tau_{\fh_0^-}\in ds, B(\tau_{\fh_0^-})\in db)$ and $p_+(z,ds,db) = p_z(\tau_{\fh_0^+}\in ds, B(\tau_{\fh_0^+})\in db)$ are hitting measures of the hypograph of $\fh_0^-(x)\coloneqq\fh_0(-x)$ and $\fh_0^+(x)\coloneqq\fh_0(x)$, $x\ge 0$, by a Brownian motion $B$ with diffusion coefficient $2$, and where $\fT_{t,x}(z_1,z_2)=\fT_{t,x}(z_1-z_2)$ with
\begin{equation}\label{eq:fTdef}
\fT_{t,x}(z)= t^{-1/3} e^{\frac{2x^3}{3t^2}-\frac{zx}{t}}\tsm\Ai(-t^{-1/3}z+t^{-4/3}x^2).
\end{equation}
This is real analytic and one checks directly (using $\Ai''(z)=z\tsm\Ai(z)$) that 
\begin{equation}\label{eq:partial-fT}
\partial_t\fT_{t,x}(z) = \tfrac13 \partial_z^3\fT_{t,x}(z),\qquad \partial_x\fT_{t,x}(z) = \partial_z^2\fT_{t,x}(z).
\end{equation}
Now fix $(x_1,\dotsc,x_n),(r_1,\dotsc,r_n)\in\rr^n$, introduce auxiliary variables $x,r\in\rr$, and let
\begin{equation}\label{eq:Khat}
\wh\fK(t,x,r)=\fK(t,x_1+x,\dotsc,x_n+x,r_1+r,\dotsc, r_n+r,\fh_0,\cdot,\cdot).
\end{equation}
The differential relations in \eqref{eq:partial-fT} together with the fact that $r$ enters the kernel just as a shift give (1), (2), (3) of Thm. \ref{thm3} for $\wh\fK$.
On the other hand, it is shown in \cite[Appx. A.1]{fixedpt} that $\|\wh\fK\|_1\longrightarrow0$ as $r\to\infty$ and that the integrals converge absolutely and uniformly in compact balls around $t>0$, $x\in\rr^n$ and $r\in\rr^n$, and thus $\wh\fK$ is real analytic in those variables as well.
Thus we have

\begin{thm}
The kernel $\wh\fK$ in \eqref{eq:Khat} satisfies the conditions of Thm. \ref{thm3}.
\end{thm}

Translating the partial derivatives in $r$ and $x$ in \eqref{eq:matKPag} back to the $r_i$'s and $x_i$'s gives \eqref{eq:matKP} for the kernel $\fK$.

There is an informal way \cite{fixedpt} to write the KPZ fixed point kernel which makes (2) and (3) of Thm. \ref{thm3} completely apparent (equation (1) comes simply again from the way the $r_i$'s come in the formula).
For $\fh\in\UC$ and $\ell_1<\ell_2$ let 
\begin{align}
{\bf P}_{\ell_1,\ell_2}^{\uptext{No hit}\,\fh}(u_1,u_2)\tts \d u_2&=\pp_{\fB(\ell_1)=u_1}\!\left(\fB(y)>\fh(y)\text{ on }[\ell_1,\ell_2],\,\fB(\ell_2)\in\d u_2\right),\\
{\bf P}_{\ell_1,\ell_2}^{\uptext{Hit}\,\fh}&=e^{(\ell_2-\ell_1)\p^2}- {\bf P}_{\ell_1,\ell_2}^{\uptext{No hit}\,\fh},
\end{align}
where $\fB$ is again a Brownian motion with diffusion coefficient $2$ and $\p$ denotes the derivativ operator.
The \emph{Brownian scattering transform of $\fh$} is the formal object
\begin{equation}\label{eq:brScatt}
\fK^{\hypo(\fh)}
 =\lim_{\ell_1\to -\infty\atop \ell_2\to \infty}e^{\ell_1\partial^2}{\bf P}_{\ell_1,\ell_2}^{\uptext{Hit}\,\fh}e^{-\ell_2\partial^2}
 =\fI-\lim_{\ell_1\to -\infty\atop \ell_2\to \infty}  e^{\ell_1\partial^2}{\bf P}_{\ell_1,\ell_2}^{\uptext{No hit}\,\fh}e^{-\ell_2\partial^2},
\end{equation}
where ${\bf P}_{\ell_1,\ell_2}^{\uptext{Hit/No hit}\,\fh}$ are thought of as operators with the given integral kernels. 
This doesn't make sense since the backward heat operator is asked to act on non-analytic functions. 
However, in the KPZ fixed point formula $\fK^{\hypo(\fh)}$ is never actually used by itself, but only after conjugation by the \emph{Airy unitary group},
\begin{equation}\label{eq:airyop}
 \fU_t=e^{ -\frac13t \partial^3}, 
\end{equation}
with $t\neq 0$ and, as we explain next, in our setting the conjugated kernel is in fact well defined.

For $a\in\rr$ define
\begin{equation}\label{eq:Pdef}
\P_a(x)=\uno{x>a}\qqand\bar\P_a(x)=\uno{x\leq a},
\end{equation}
which we also regard as multiplication operators acting on $L^2(\rr)$.
For $t>0$ the Airy semigroup acts by convolution with Airy functions, i.e. $\fU_t$ has integral kernel $\fU_t(x,y)=t^{-1/3}\Ai(t^{-1/3}(x-y))$; its inverse $\fU_t^{-1}=\fU_{-t}$ equals $\fU_t^*$.
The Airy functions are not themselves in $L^2(\rr)$; however, for $t>0$ and $r>-\infty$, $ \fU_t^{-1}\P_{r}$ maps $L^2(\rr)$ into the domain of $e^{x \partial^2}$ for \emph{any} $x\in\rr$.
So for $t>0$ and $r>-\infty$, we \emph{define} on $L^2([r,\infty))$
\begin{equation}
\fK^{\hypo(\fh)}_t=  \lim_{\ell_1\to -\infty\atop \ell_2\to \infty}
\fU_t  e^{\ell_1\partial^2}{\bf P}_{\ell_1,\ell_2}^{\uptext{Hit}\,\fh}e^{-\ell_2\partial^2}\fU_t^{-1}.\label{def:scattering}
\end{equation}
For any $t>0$ and $r>-\infty$ the limit on the right hand side of \eqref{def:scattering} exists in trace class on $L^2([r,\infty))$ \cite{flat,fixedpt}, and defines the left hand side as a trace class operator in this space.  
The limit is given by the right hand side of \eqref{eq:trcldecomp2} with $x_i=x_j=x$ and $r_i=r_j=0$; the convergence was proved first in \cite{flat}.
It satisfies the semigroup property  
\begin{equation}
\fU_s\fK^{\hypo(\fh)}_t\fU_s^{-1}=\fK^{\hypo(\fh)}_{\ft+\fs}
\end{equation}
so we can write (at least informally)
\begin{equation}\label{sg}
\fK^{\hypo(\fh)}_t=\fU_t\fK^{\hypo(\fh)}\fU_t^{-1}.
\end{equation}
Note that we avoid the problem of domains by not defining the left hand side of \eqref{def:scattering} as a product of three operators, but just as one operator with the semigroup property.  In this sense the Brownian scattering operator is the germ of the semigroup.
Alternatively one can think of the Brownian scattering operator as the entire semigroup $(\fK^{\hypo(\fh)}_t)_{t>0}$ from \eqref{sg}.
The fact that \eqref{eq:brScatt} is formal is important, though.
We will see in \eqref{eq:K0} that \emph{the limit of \eqref{sg} as $t\searrow 0$ is not $\fK^{\hypo(\fh)}$}.

From $\fK^{\hypo(\fh)}$ we build an \emph{extended Brownian scattering operator} acting on $L^2(\{x_1,\cdots,x_m\}\times\rr)$, 
\begin{equation}\label{eq:scatt-ext}
\fK^{\hypo(\fh)}_{\uptext{ext}}(x_i,\cdot;x_j,\cdot)=-e^{(x_j-x_i)\p^2}\uno{x_i<x_j}+e^{-x_i\p^2}\fK^{\hypo(\fh)} e^{x_j\p^2},
\end{equation}
with the analogous caveat that in order to make sense  of it each of the above $(x_i,x_j)$ entries should be conjugated by $\fU_t$ and surrounded by $\P_{r_i}$ on the left and $\P_{r_j}$ on the right.
Then the \emph{KPZ fixed point kernel} appearing on the right hand side of \eqref{eq:F1a} and given in \eqref{eq:trcldecomp2} can be written as
\begin{equation}
\fK_{ij}(u_1,u_2)=\fK^{\hypo(\fh)}_{t,\uptext{ext}}(x_i,u_1+r_i;x_j,u_2+r_j)\quad\uptext{with}\quad \fK^{\hypo(\fh_0)}_{t,\uptext{ext}}=\fU_t\fK^{\hypo(\fh_0)}_{\uptext{ext}}\fU_t^{-1},\label{eq:KKscatt}
\end{equation}
where the Airy operators act on the left and right of each entry of $\fK^{\hypo(\fh_0)}_{\uptext{ext}}$.
From this expression one sees that $\p_t\fK_{ij}=-\frac13(\p^3\fK_{ij}-\fK_{ij}\p^3)$ and $\p_{x_\ell}\fK_{ij}=-\frac12(\p^2\fK_{ij}\uno{i=\ell}-\fK^2_{ij}\p^2\uno{j=\ell})$ which, after integration by parts and summing over the $x_\ell$'s as above, become (2) and (3) in Thm. \ref{thm3}.

\begin{rem}
For compactly supported initial data, which in our context means that $\fh(y)=-\infty$ for $y$ outside some compact interval, the definition of the Brownian scattering operator in \eqref{def:scattering} becomes explicit, without any limits.
In fact, consider $\fh$ with support on some interval $[-a,a]$ and let $\ell_1<-a$ and $\ell_2>a$.
Clearly ${\bf P}_{\ell_1,\ell_2}^{\uptext{Hit}\,\fh}=e^{-(a+\ell_1)\p^2}{\bf P}_{-a,a}^{\uptext{Hit}\,\fh}e^{(\ell_2-a)\p^2}$, so $\fU_te^{\ell_1\partial^2}{\bf P}_{\ell_1,\ell_2}^{\uptext{Hit}\,\fh}e^{-\ell_2\partial^2}\fU_t^{-1}=\fU_te^{-a\p^2}{\bf P}_{-a,a}^{\uptext{Hit}\,\fh}e^{-a\p^2}\fU_t^{-1}$.
Since $\fU_te^{x\p^2}=\fT_{t,x}$ and $\fU_t^{-1}e^{x\p^2}=\fT_{-t,x}$ (see \eqref{eq:partial-fT}), this shows that
\begin{equation}
\fK^{\hypo(\fh)}_t=\fT_{t,-a}{\bf P}_{-a,a}^{\uptext{Hit}\,\fh}\fT_{-t,-a}.
\end{equation}
The right hand side coincides with the right hand side of \eqref{eq:trcldecomp2} with $i=j$ and $x_i=r_i=0$ (see \cite[Sec. 4.1]{fixedpt}).
\end{rem}

\subsection{Initial data} 

The natural class of initial data for our problem (the ``one dimensional substrate'') corresponds to functions in the space\footnote{With some work the growth condition on the initial data can be relaxed to  $\fh(x) \le A x^2 + B$ for the problem up to a finite time $t=t(A)$.}
\begin{multline}
\UC=\big\{\fh\!:\rr\longrightarrow[-\infty,\infty)\!:\fh\text{ is upper semicontinuous,}\\  
\fh(x)\leq A+B|x|\text{ for some $A,B>0$ }\text{and }\fh\not\equiv-\infty\big\}.\label{eq:UC}
\end{multline}
A function is upper semicontinuous if and only if its \emph{hypograph} $\hypo(\fh) = \{(\fx,y): y\le \fh(\fx)\}$ is closed in $[-\infty,\infty)\times \rr$.
We endow $[-\infty,\infty)$ with the distance $d_{[-\infty,\infty)}(y_1,y_2) = |e^{y_1} - e^{y_2}|$, and 
use the topology of \emph{local Hausdorff convergence}, which means Hausdorff convergence of the restrictions to $-L\le x\le L$ of 
$\hypo(\fh_n)$ to $\hypo(\fh)$ for each $L>0$. 

\begin{ex}{\bf (Finite collection of narrow wedges)}\label{fcnw}
\enspace
Let $a_1<a_2<\dotsm<a_k$, $b_1,b_2,\dotsc,b_k\in\rr$.
Then $\fh=\mathfrak{d}_{\vec a}^{\vec b}$ is in $\UC$, with $\mathfrak{d}_{\vec a}^{\vec b}$ defined by 
\[\mathfrak{d}_{\vec a}^{\vec b}(x)=b_i\quad\uptext{if $x=a_i$ for some $i$},\qquad\mathfrak{d}_{\vec a}^{\vec b}(x)=-\infty\quad\uptext{otherwise}.\]
\end{ex}

The initial data for the one point case \eqref{eq:KP-II} is the \emph{``escarpment''}
\begin{equation}\label{eq:escarpment}
\phi(0,x,r) = 0\quad\uptext{for }r\ge \fh_0(x),\qquad\phi(0,x,r)=-\infty\quad\uptext{for }r<\fh_0(x).
\end{equation}
These are unusual and do not fit into any well-posedness schemes known for the KP equation\footnote{\cite{kenigPonceVega} consider as initial data for KdV an odd polynomial with positive leading coefficient, which is somewhat in the same spirit.}.  Although the $\infty$ looks formal, we believe the solutions to the equations with such initial data  are well posed, but we leave the proofs for future work.  They also appear not to develop solitons.
Since $F$ is given by a Fredholm determinant, these  initial conditions represent an entirely new class of integrable initial data for KP. 

The initial data for the matrix KP equation \eqref{eq:matKP} is formally
\begin{equation}\label{eq:matKPini}
Q_{i,j}(0,x_1,\dotsc,x_n,r_1,\dotsc,r_n)=
\begin{dcases*}
-{\bf P}^{{\geq\fh_0,\leq -\mathfrak{d}_{\vec x}^{-\vec r}}}_{x_i,x_j}(r_i,r_j) & if $i<j$,\\
\infty & if $i=j$ and $r_i<\fh_0(x_i)$,\\
0 & otherwise.
\end{dcases*}
\end{equation}
where
\begin{align}\label{eq:matrix0}
&{\bf P}^{{\geq\fh_0,\leq -\mathfrak{d}_{\vec x}^{-\vec r}}}_{x_i,x_j}\!(r_i,r_j)dr_j\\
&\quad=\pp_{\fB(x_i)=r_i}\!\big(\fB(y)\geq \fh_0(y)~\forall\,y\in[x_i,x_j],\fB(x_n)\leq r_n\uptext{ for each }x_n\in(x_i,x_j), \fB(x_j)\in dr_j\big).
\end{align}
The probability is with respect to a Brownian motion $\fB$ with diffusivity $2$ starting at $r_i$ at time $x_i$.
This is derived in Appx. \ref{mid} for finite collections of narrow wedges.  Unlike the scalar case, one can see immediately that the initial data as written is insufficient because in the matrix product the $0$ and $\infty$ interact.
Wrongly interpreted it appears to produce anomalous solutions, so the initial data written
would have to be augmented by at least some description of the rate of convergence to $0$ and $\infty$ in the $t\downarrow 0$ limit.
We leave this also for future work.

\begin{rem}
 When $t=0$, the escarpment initial data \eqref{eq:escarpment} just means (sending $x$ to $-x$) that $L$ from \eqref{lax} is the heat operator with Dirichlet boundary data on the hypograph of $\fh_0$.  This corresponds to the Brownian scattering transform which computes transition probabilities of Brownian motions killed when passing through that hypograph.
\end{rem}

\noindent {\bf Outline.}\enspace In Sec. \ref{sec:examples} we show how the famous Painlev\'e expressions for the Tracy-Widom distributions just arise as special examples of the KP equation, and write equations for a few other special initial data, in particular the Airy processes.  In Sec.  \ref{sec3} we show that some of the explicit formulas for the one-dimensional distributions of the KPZ equation (as opposed to the KPZ fixed point) also satisfy the KP equation.
Sec. \ref{sec:proof} contains the proof of Thm. \ref{thm3}.

\section{Examples}\label{sec:examples}

\subsection{Tracy-Widom distributions}\label{sec:TW} 

A key observation is that \emph{the GUE and GOE Tracy-Widom distributions are now seen to simply arise as special similarity solutions of the KP equation \eqref{eq:KP-II}}:

\begin{ex}{\bf (Tracy-Widom GUE distribution)}\label{ex:GUE}
\enspace
Consider $\fh_0=\mathfrak{d}_0$, the \emph{narrow wedge initial condition} defined as $\mathfrak{d}_0(0)=0$ and $\mathfrak{d}_0(x)=-\infty$, $x\neq 0$.
With this choice of initial data one has $\fh(t,x)+x^2/t\stackrel{\uptext{dist}}{=}t^{1/3}\aip(t^{-2/3}x)$ where $\aip$ is the Airy$_2$ process (see Sec. \ref{sec:airy}), which is stationary (in $x$).
From this and the 1:2:3 scaling invariance of \eqref{eq:KP-II}, it is natural to look for a self-similar solution of the form
\begin{equation}
\phi^\uptext{nw} (t,x,r)=t^{-2/3} \psi^\uptext{nw} (t^{-1/3} r + t^{-4/3}x^2).
\end{equation}
This turns \eqref{eq:KP-II} into
\begin{equation}\label{eq:KP-IIA}
(\psi^\uptext{nw})'''+12\psi^\uptext{nw}(\psi^\uptext{nw})'-4r(\psi^\uptext{nw})'-2\psi^\uptext{nw}=0.
\end{equation}
The transformation $\psi^\uptext{nw}=-\q^2$ takes \eqref{eq:KP-IIA} into
Painlev\'e II:
\begin{equation}\label{eq:p2}
\q''=r\q+2\q^3.
\end{equation}
As $r\to -\infty$ the solution is approximately $\phi^\uptext{nw}(t,x,r)\sim -(\tfrac{r}{2t} +\tfrac{x^2}{2t^2})$, picking out the Hastings-McLeod solution $\q(r)\sim -\Ai(r)$ as $r\to \infty$.
Thus we recover
\begin{equation}
F(t,x,r)= F_\uptext{GUE}(t^{-1/3} r+t^{-4/3}x^2)
\end{equation}
where $F_\uptext{GUE}$ is the GUE Tracy-Widom distribution \cite{tracyWidom}, usually written in the equivalent form
\begin{equation}
F_\uptext{GUE}(s) = \exp\!\left\{ -\int_s^\infty\!\d u\, (u-s)\q^2(u)\right\}.
\end{equation}
\end{ex}

\begin{ex}{\bf (Tracy-Widom GOE distribution)}\label{ex:GOE}
\enspace
If $\fh_0(x)\equiv 0$, the \emph{flat initial condition}, there is no $x$ dependence and \eqref{eq:KP-II} reduces to KdV. 
Now we look for a self-similar solution of \eqref{eq:KdV} the form
\begin{equation}\label{flatform}
\phi^\uptext{fl}(t,r) = (t/4)^{-2/3} \psi^\uptext{fl}((t/4)^{-1/3} r)
\end{equation}
(the extra factor of $1/4$ is just to coordinate conventions with random matrix theory), obtaining the ordinary differential equation
\begin{equation}
(\psi^\uptext{fl})''' +12(\psi^\uptext{fl})'\psi^\uptext{fl}-r(\psi^\uptext{fl})' -2\psi^\uptext{fl} =0.
\end{equation}  
Miura's transform 
\begin{equation}
\psi^\uptext{fl} = \tfrac12 ( \q' -\q^2)
\end{equation}  
brings this to Painlev\'e II  \eqref{eq:p2}, with the same behavior as $r\to-\infty$.
So we recover
\begin{equation}
F(t,x,r)= F_\uptext{GOE}(4^{1/3}t^{-1/3} r)
\end{equation}
where $F_\uptext{GOE}$ is the GOE Tracy-Widom distribution \cite{tracyWidom2}, usually written in the equivalent form
\begin{equation}
F_\uptext{GOE}(r) =  \exp\!\left\{ -\frac12\int_r^\infty\!\d u\,\q(u)\right\}F_\uptext{GUE}(r)^{1/2}.
\end{equation}
\end{ex}

These two examples also have the following interpretation. 
Let $\lambda^\uptext{max,GUE}_N$ and $\lambda^\uptext{max,GOE}_N$ be the largest eigenvalues of $N\times N$ matrices chosen from the Gaussian Unitary and Gaussian Orthogonal Ensembles multiplied by $\sqrt{N}$, so that $\lambda^\uptext{max,GUE}_N\sim 2N + N^{1/3}\zeta_\uptext{GUE}$ and  $\lambda^\uptext{max,GOE}_N\sim 2N + N^{1/3}\zeta_\uptext{GOE}$ with $\zeta_\uptext{GUE}$ and $\zeta_\uptext{GOE}$ the standard Tracy-Widom GUE and GOE random variables.
Let 
\begin{align}
F_1(t,r)& = {\textstyle\lim_{N\to \infty}} \pp(N^{-1/3}( \lambda^\uptext{max,GOE}_{Nt}-2Nt) \le 4^{1/3}r)=F_\uptext{GOE}(4^{1/3}t^{-1/3}r),\\
F_2(t,x,r) &= {\textstyle\lim_{N\to \infty}} \pp(N^{-1/3}( \lambda^\uptext{max,GUE}_{Nt}-2 Nt) \le r+x^2/t)=F_\uptext{GUE}(t^{-1/3} r+t^{-4/3}x^2).
\end{align}
As we have seen, $\p_r^2\log F_1$ and $\p_r^2\log F_2$ satisfy the KP equation \eqref{eq:KP-II}.
In the former case, there is no dependence on $x$ and KP reduces to KdV \eqref{eq:KdV}.

\subsection{PDEs for other initial data}

Another question is whether there are analogues of Painlev\'e II for other self-similar solutions.
It is natural to observe $\phi$ in the frame of the inviscid solution $\tfrac14(\p_x\bar\fh)^2-\p_t\bar\fh=0$ of Burgers' equation,
\begin{equation}
\phi(t,x,r) \coloneqq \bar\phi(t,x,r-\bar\fh(t,x));
\end{equation}
one obtains
\begin{equation}\label{stan}\partial_t  \bar\phi +\bar\phi\ts\partial_r  \bar\phi + \tfrac1{12}\partial_r^3  \bar\phi + \tfrac14\partial_r^{-1} \partial_x^2  \bar\phi -\tfrac12\p_x\bar\fh\ts\p_x\bar\phi+ V \bar\phi= 0\end{equation}
with $V=-\tfrac{1}4\p^2_x\bar\fh$ and with initial data $\bar\phi(0,x,r) = 0$ for $r\ge 0$ and $-\infty$ for $r<0$. 

In order to get a solution for the rescaled spatial process, write
\begin{equation}
\phi(t,x,r) = t^{-2/3}\psi(t,t^{-2/3}x,t^{-1/3}(r- \bar\fh(t,x))).
\end{equation}
Then
\[t\tts\p_t\psi-\tfrac23\psi-\tfrac13r\p_r\psi+\tfrac1{12}\p_r^3\psi+\psi\p_r\psi+\tfrac14\p_r^{-1}\p_x^2\psi -\tfrac{1}4\p_x^2\tilde\fh\tts\psi -(\tfrac23x +\tfrac12\p_x\tilde\fh)\p_x\psi =0\]
with $\tilde\fh(t,x)=t^{-1/3}\bar\fh(t,t^{2/3}x)$.

\begin{ex}{\bf (Half-flat initial data)}
\enspace Consider $\fh_0(x)=0$, $x\le0$ and $\fh_0(x)=-\infty$, $x>0$.
Now $\bar\fh(t,x) = -x^2/t\tts\uno{x\geq0}$.
There is dependence on $x$, though not on $t$.
This gives rise to a partial differential equation for $\psi(x,r)$,
\[-\tfrac13r\p_r\psi+\tfrac1{12}\p_r^3\psi-(\tfrac16\uno{x\ge0}+\tfrac23\uno{x<0})\psi+\psi\p_r\psi+\tfrac14\p_r^{-1}\p_x^2\psi +(\tfrac13\uno{x\ge 0}-\tfrac23\uno{x<0})x\p_x\psi=0.\]
\end{ex}

\begin{rem}{\bf (Lower tail heuristics)}
\enspace 
Typically the equation is controlled on large scales by the equation with the third derivative dropped, and the Burgers' equation makes sense for such wedge type initial data.
Let $\bar\phi$ be as in \eqref{stan} and
\begin{equation}
\check\phi = \bar\phi -\eta  
\end{equation}
where $\eta=\tfrac{c }{t} r\uno{r<0}$.
Then $\bar\phi(0,x,r)=0$ and
\begin{equation}\label{eq:KP-IIb}
\partial_t \check\phi + \tfrac12\partial_r \check\phi^2+ \tfrac1{12}\partial_r^3\check\phi+ \tfrac14\partial_r^{-1} \partial_x^2\check\phi + \partial_r(\eta\check\phi)
-\tfrac12 \partial_x \bar\fh\p_x\check\phi-\tfrac14 \partial^2_x \bar\fh~ \check\phi +(\tfrac{c-1}{t} -\tfrac14 \partial^2_x \bar\fh) \eta - \tfrac{c }{12t} \delta'_0(r)
= 0.
\end{equation}
One hopes to set things up so that $\check\phi$ has good decay at $\pm\infty$.
Consider our two basic examples.
In the flat case $\bar\fh\equiv0$, and we want to take $c=1$ to make the second last term drop out, which will lead to the conclusion that $\phi(t,r)\sim r/t$, or $\log F(t,r)\sim - \tfrac1{6t} r^3$ as $r\to-\infty$, recovering the Tracy-Widom GOE lower tail. 
In the narrow wedge case $\bar\fh = -\tfrac{x^2}{t}$ so $\tfrac14 \partial^2_x \bar\fh=-\tfrac{1}{2t}$ and we want to take $c=1/2$ leading to $\log F(t,x,r)\sim  -\tfrac1{12 t} (r+ \tfrac{x^2}{t})^3 $, recovering the Tracy-Widom GUE lower tail.  

\noindent The conclusion is that \emph{the lower tail of the distributions can be seen directly from the ``Burgers'' part of the KP equation, which dominates in that region}.

\noindent Note that this is only intended to be a quick heuristic.  If one were to try to make it rigorous, say, for the flat and narrow wedge case, one would be led through self-similar solutions back to the Painlev\'e equation and the standard tools there to get the left tails (see e.g. \cite{baikBuckDiF}).
Such Riemann-Hilbert methods give much finer information, for example, subleading terms.
It is an interesting question whether such results can be obtained from the PDE \eqref{stan} so as to get finer information on the left tail for more general initial data.
\end{rem}

\subsection{Airy process} \label{sec:airy}

The \emph{Airy process} $\aip(x)$ is defined as 
\begin{equation} 
\aip(x) \coloneqq \fh(1,x;\mathfrak{d}_0)+x^2
\end{equation}
where $\fh(t,x;\mathfrak{d}_0)$ is the KPZ fixed point starting from a narrow wedge $\mathfrak{d}_0$ at the origin.
The Airy process is stationary and the one point distribution is the Tracy-Widom GUE distribution.

K. Johansson famously asked whether there is an equation for the multipoint distribution.
PDEs were given by \cite{adlerVanMoPDEs,tracyWidom-AiryDiffEqns}.
Starting with narrow wedge initial data the matrix KP equation \eqref{eq:matKP} gives us another PDE.
This equation lives at the same level as the system of PDEs in \cite{tracyWidom-AiryDiffEqns}; it could well be that they are equivalent, but we have not succeeded in confirming it yet, and leave it for future work.

Note the KP-II equation is written in variables $t$, $r=r_1+\dotsc+r_n$ and $x=x_1+\dotsc+x_n$, the rest of the variables entering only from the boundary condition.
Exploiting skew time reversal invariance we find an extra symmetry in the narrow wedge case.
Let $\mathfrak{d}^{\vec r}_{\vec x}$ denote a multiple narrow wedge (see Ex. \ref{fcnw}).
For $\vec x\in\rr^m$ and $z\in\rr$ write $\vec x+z=(x_1+z,\dotsc,x_m+z)$.
By skew time reversibility and translation and affine invariance of the KPZ fixed point \cite[Thm. 4.5]{fixedpt},
\begin{equation}\label{eq:sktime-nw}
\begin{split}
F(t,\vec x+z,\vec r+a)&=\pp\big(\fh(t,x_i+z;\mathfrak{d}_0)\leq r_i+a\text{ for all }i\big)
=\pp\big(\fh(t,\cdot;\mathfrak{d}_0)\leq-\mathfrak{d}^{-\vec r-a}_{\vec x+z}\big)\\
&=\pp\big(\fh(t,\cdot;\mathfrak{d}^{-\vec r-a}_{\vec x+z})\leq-\mathfrak{d}_0\big)
=\pp\big(\fh(t,0;\mathfrak{d}^{-\vec r}_{\vec x+z})\leq a\big)\\
&=\pp\big(\fh(t,-z;\mathfrak{d}^{-\vec r}_{\vec x})\leq a\big).
\end{split}\end{equation}
Now the right hand side is just the one point distribution at $-z$ with a given, fixed initial condition.
So if we let $G(t,z,a)=F(t,\vec x+z,\vec r+a)$ we see that \emph{$\p_a^2\log G$ satisfies \eqref{eq:KP-II} in $(t,z,a)$} (in terms of \eqref{eq:escarpment}, the initial data is $G(0,z,a)=-\infty$ if $z=-x_i$ for some $i$ and $a<-r_i$ and $G(0,z,a)=0$ otherwise).
But $\p_zG(t,z,a)=\cD_xF(t,\vec x+z,\vec r+a)$ and similarly $\p_aG(t,z,a)=\cD_rF(t,\vec x+z,\vec r+a)$.
So setting now $a=z=0$ we deduce that $\phi=\cD_r^2\log F$ satisfies
\begin{equation}
\partial_t \phi + \tfrac12 \cD_r \phi^2 + \tfrac1{12}\cD_r^3 \phi + \tfrac14\cD_r^{-1} \cD_x^2 \phi = 0\label{eq:KP-II-Airy}
\end{equation}
(the initial data is now similarly $\phi(0,\vec x,\vec r)=-\infty$ if $x_i=0$ for some $i$ and $r_i<0$ and $\phi(0,\vec x,\vec r)=0$ otherwise).

Note that from Thm. \ref{thm:2}, we know that $\cD_r^2\log F=\tr q$ and $\p_t\tr q+\tr(q\cD_rq)+\tfrac1{12}\cD_r^3\tr q+\tfrac14\cD_x^2\tr Q=0$ (using $\tr(AB)=\tr(BA)$).
But the above argument implies then that $\tr q$ \emph{itself solves KP} \eqref{eq:KP-II-Airy}, and
as a consequence, we deduce in the narrow wedge case that $\tr(q\cD_rq)=\tr(q)\tr(\cD_rq)$.
This can also be proved directly using the fact that, in this case, $(\p_u+\p_v)\fK^{\hypo(\fh_0)}_t(u,v)$ is a rank one kernel, which implies that $q$ is a rank one matrix.

An alternative derivation of \eqref{eq:KP-II-Airy} using the path integral formula for the KPZ fixed point can be found in Appendix \ref{pif}.

Next we use the 1:2:3 scaling invariance of the KPZ fixed point and stationarity of the Airy process, as in Ex. \ref{ex:GUE}. 
Let $H$, $\Psi$, $\psi$ denote the $(m+1)$-point distribution function of the Airy process and its logarithmic derivatives,
\begin{gather}
H(y_1,\dotsc,y_{m},r_0,\dotsc,r_m) = \pp( \aip(0)\le r_0,\aip(y_1)\le r_1,\dotsc,\aip(y_m)\leq r_m),\label{m-pt}
\\ \Psi= \cD_r \log H,\qquad \psi=\cD_r \Psi.\label{m-pt-der}
\end{gather}
Since $F(t,x_0,\dotsc,x_m,r_0,\dotsc,r_m)=H(t^{-2/3}\bar x_1,\dotsc,t^{-2/3}\bar x_m, t^{-1/3} (r_0 + \tfrac{x^2_0}{t}), \dotsc,t^{-1/3} (r_m +  \tfrac{x^2_m}{t}))$ with $\bar x_i=x_i-x_0$, \eqref{eq:KP-II-Airy} leads to

\begin{thm}{\bf (Airy$_2$ process multipoint function)}
\enspace The logarithmic derivatives \eqref{m-pt-der} of the multipoint function of the Airy process satisfy the PDE
\begin{multline}\label{eqfora-m}
\textstyle-\tfrac1{6}\psi + \tfrac12\cD_r\psi^2 + \tfrac1{12} \cD_r^3 \psi -\tfrac13\sum_{a=0}^mr_a\p_{r_a}\psi-\tfrac23\sum_{a=1}^my_a\p_{y_a}\psi\\
\textstyle-\sum_{a=1}^m\sum_{b=0,\,b\neq a}^my_a^2\p_{r_a}\tsm\p_{r_b}\Psi+\sum_{a,b=1,\,a\neq b}^my_ay_b\tts\p_{r_a}\tsm\p_{r_b}\Psi =0
\end{multline}
with boundary condition $\psi(0,\dotsc,0,r_1,\dotsc,r_m) =\psi^\uptext{nw}( \min(r_1,\dotsc,r_m))$ from  \eqref{eq:KP-IIA}.
\end{thm}

In particular, for the logarithmic derivative of the Airy process two-point function $\Psi(y,r_1,r_2)=\cD_r\log \pp(\aip(0)\leq r_1,\,\aip(y)\leq r_2)$ we get
\begin{equation}\label{eqfora-2}
-\tfrac1{6}\psi+\tfrac12\cD_r\psi^2+\tfrac1{12}\cD_r^3\psi-\tfrac13(r_1\p_{r_1}\psi+r_2\p_{r_2}\psi)-\tfrac23y\tts\p_y\psi-y^2\p_{r_1}\p_{r_2}\Psi =0.
\end{equation}
In our notation, the equation derived in \cite{adlerVanMoPDEs} reads
\begin{equation}\label{avm}
(r_2-r_1)\p_{r_1}\p_{r_2}\psi+y\tts\p_y(\p_{r_1}-\p_{r_2})\psi+y^2(\p_{r_1}-\p_{r_2})\p_{r_1}\p_{r_2}\Psi
+\p_{r_1}\!\Psi\ts\p_{r_2}\tsm\cD_r\psi-\p_{r_2}\!\Psi\ts\p_{r_1}\tsm\cD_r\psi=0.
\end{equation}
Although \eqref{eqfora-2} and \eqref{avm} have similarities, they do not appear to 
be equivalent.  We have not succeeded in reconciling them and leave this for future work.
Of course the Airy process has many symmetries and it is plausible that the equations are 
not equivalent yet both hold.

\begin{ex}{\bf (Airy$_1$ process multipoint function)}\label{ex:a1}
\enspace
Consider now the KPZ fixed point with flat initial data, $\fh_0\equiv 0$.
The multipoint function satisfies
\[F(t,\vec x,\vec r)=\pp\big(\aipo(0)\leq t^{-1/3}r_1,\,\aipo(t^{-2/3}\bar x_2)\leq t^{-1/3}r_2,\dotsc,\aipo(t^{-2/3}\bar x_m)\leq t^{-1/3}r_m\big)\]
with $\bar x_i=x_i-x_1$, where $\aipo(x)\coloneqq\fh(1,x;0)$ is the Airy$_1$ process and we have used the fact that it is stationary.
We know that the left hand side equals $\tr Q$ where the $m\times m$ matrix $Q$ and its derivative $q$ solve \eqref{eq:matKP}.
But the right hand side only depends on the differences $x_i-x_1$, so $Q$ actually solves the (integrated) matrix KdV equation
\[\p_tQ+\tfrac12(\cD_rQ)^2+\tfrac1{12}\cD_r^3Q=0.\]
Take now $m=2$ for simplicity.
The (formal) initial data $Q(0,x_1,x_2,r_1,r_2)$, \eqref{eq:matKPini}, is also invariant under $r_i,x_i\mapsto t^{-1/3}r_i,t^{-2/3}x_i$ followed by $Q\mapsto t^{-1/3}Q$, and only depends on $|x_2-x_1|$ (in fact $Q_{12}(0,r_1,r_2,x_1,x_2) = \uno{r_1\ge0,r_2\ge0}\frac{1}{\sqrt{4\pi|x_2-x_1|}}\big(e^{-(r_1+r_2)^2/4|x_2-x_1|}-e^{-(r_1-r_2)^2/4|x_2-x_1|}\big)$ by the reflection principle).
In this case we look for a matrix solution of \eqref{eq:matKP} with the scaling
\[Q(t,r_1,r_2,x_1,x_2)= t^{-1/3}\bar{Q}(t^{-2/3}(x_2-x_1),t^{-1/3}r_1,t^{-1/3}r_2).\]
The conclusion is that
\begin{equation}
\cD_r\log\pp(\aipo(0)\leq r_1,\,\aipo(y)\leq r_2)=\tr\!\big(\bar Q(y,r_1,r_2)\big)\label{eq:2ptAiry1Q}
\end{equation}
with $\bar Q(y,r_1,r_2)$ solving the matrix PDE
\[-\tfrac13\bar Q+\tfrac12(\cD_r\bar Q)^2+\tfrac1{12}\cD_r^3\bar Q-\tfrac13(r_1\p_{r_1}\bar Q+r_2\p_{r_2}\bar Q)-\tfrac23y\p_y\bar Q=0.\]
It does not seem to be possible to turn this $4\times4$ system into a closed equation for the left hand side of \eqref{eq:2ptAiry1Q}; the fact that this works for the Airy$_2$ process is very particular to narrow wedge initial data (and follows from skew time reversibility as used in \eqref{eq:sktime-nw}).
\end{ex}

\settocdepth{section}
\hypersetup{bookmarksdepth=1}

\section{KP-II in special solutions of the KPZ equation} \label{sec3}

Thm. \ref{thm3} also implies that \emph{some of the special explicit solutions for one dimensional distributions of the Kardar-Parisi-Zhang equation \eqref{KPZ} satisfy the (scalar) KP-II equation \eqref{eq:KP-II} as well}.
At this point we do not know if this is part of a more general fact, or if KP-II only holds in these special cases because of some symmetry.
All we have is examples.

\subsection*{Narrow wedge solution of the KPZ equation}\label{ex:nw}
Let $\kh_\uptext{nw}$ be the narrow wedge solution of \eqref{KPZ} with $\lambda=\nu=\tfrac14$ and $\sigma=1$.
In other words, $\kh_\uptext{nw}=\log Z$ where $Z$ is the fundamental solution of the stochastic heat equation (SHE) with multiplicative noise
\begin{equation}\label{eq:SHE}
\p_t Z = \tfrac14 \p^2_xZ + \xi Z,\qquad Z(0,x)=\delta_0(x).
\end{equation}
The KPZ generating function is
\begin{equation}\label{eq:gumbel}
G_\uptext{nw}(t,x,r) = \ee\!\left[\exp\!\left\{-e^{\kh_\uptext{nw}(t,x)+\frac{t}{12}-r}\right\}\right].
\end{equation}
The distribution of $\kh_\uptext{nw}(t,x)$ was computed in 2010 in \cite{acq,sasamSpohn,dotsenkoNW,calabreseLeDousallNW}, with the result that\footnote{See \cite[Thm. 1.10]{borodinCorwinFerrari}, which computes the generating function directly. In comparing with that formula, we are changing variables $(t,x)\longmapsto(2t,2x)$ to match the two different scaling conventions for \eqref{eq:SHE}, and using the fact that $\kh_\uptext{nw}(t,x)+x^2/t$ is stationary.\label{ft:sc}} $G_\uptext{nw}(t,x,r)=\det(\fI - \fK)_{L^2(\rr^+)}$
with
\[\fK(u,v)=\int_{-\infty}^\infty dy\,t^{-2/3}\frac{1}{1+e^{y}}\Ai(t^{-1/3}(u+r-y)+t^{-4/3}x^2)\Ai(t^{-1/3}(v+r-y)+t^{-4/3}x^2).\]
If we conjugate the operator by multiplying the kernel by $e^{(v-u)x/t}$ the value of the Fredholm determinant doesn't change, i.e. we also have 
\[G_\uptext{nw}(t,x,r)=\det(\fI - \wt\fK)_{L^2(\rr^+)}\quad\uptext{with}\quad\wt\fK(u,v)=e^{(v-u)x/t}\fK(u,v).\]
One checks directly that $\wt\fK$ satisfies the differential relations (1)--(3) from Thm. \ref{thm3} (this can also be seen by writing $\wt\fK(u,v) = \fU_te^{-x\p^2}\tts\fM\tts e^{x\p^2}\fU_t^{-1}(u+r,v+r)$ with $\fM$ the multiplication operator $\fM f(u) = (1+e^{u})^{-1}f(u)$ and $\fU_t$ the Airy unitary operator defined in \eqref{eq:airyop}).
It follows that:

\begin{thm}
 $\phi_\uptext{nw}\coloneqq\p_r^2\log G_\uptext{nw}$ satisfies the KP-II equation \eqref{eq:KP-II}.
\end{thm}

The initial condition is  $\lim_{t\searrow0}\phi_\uptext{nw}(t,x,r-\tfrac{x^2}{t} -\log \sqrt{\pi t} ) =  -e^{-r}$.  This suggests defining the $x$ independent, shifted variable $\hat \phi_\uptext{nw}(t,r )= \phi_\uptext{nw}(t,x,r-\tfrac{x^2}{t} -\log \sqrt{\pi t} )$ which now satisfies the \emph{cylindrical KdV equation},
\begin{equation}
\partial_t \hat \phi_\uptext{nw} +\tfrac1{2t}\partial_r \hat \phi_\uptext{nw}+  \hat \phi_\uptext{nw}\partial_r \hat \phi_\uptext{nw}+ \tfrac1{12}\partial_r^3 \hat \phi_\uptext{nw}  +\tfrac1{2t} \hat \phi_\uptext{nw}= 0,\qquad  \hat \phi_\uptext{nw}(0,r) =-e^{-r} .
\end{equation}


\subsection*{KPZ equation with spiked/half-Brownian initial data}

Consider now the solution $\kh_{b}$ of \eqref{KPZ} with $\lambda=\nu= \tfrac14$, $\sigma=1$, and $m$-spiked initial data, where $b=(b_1,\dotsm,b_m)\in\rr^m$ are the spike parameters.
When $m=1$, this corresponds to half-Brownian initial data (more precisely, at the level of the SHE one sets $Z(0,x)=e^{B(x)+b_1x}\uno{x\geq0}$ where $B(x)$ is a Brownian motion with diffusivity 2); for the definition in the general case $m\geq1$  we refer to \cite[Defn. 1.9]{borodinCorwinFerrari}.
Define $G_b$ as in \eqref{eq:gumbel} (with $\kh_b$ in place of $\kh_\uptext{nw}$).
Then from \cite[Thm. 1.10]{borodinCorwinFerrari} we get now that $G_b(t,0,r)=\det(\fI - \fK_0)_{L^2[0,\infty)}$ with (see footnote \ref{ft:sc} again)
\[\fK_0(u,v)=\itwopii{2}\int_{\mathcal{C}_t} d\eta\int_{\mathcal{C}_t'} d\xi\,\frac{t^{-1/3}\pi}{\sin(t^{-1/3}\pi(\xi-\eta))}\frac{e^{\xi^3/3-(u+t^{-1/3}r)\xi}}{e^{\eta^3/3-(v+t^{-1/3}r)\eta}}\prod_{k=1}^m\frac{\Gamma(t^{-1/3}\eta-b_k)}{\Gamma(t^{-1/3}\xi-b_k)},\]
where $\mathcal{C}_t$ goes from $-\frac14t^{1/3}-\I\infty$ to $-\tfrac14t^{1/3}+\I\infty$ crossing the real axis to the right of $t^{1/3}b_1,\dotsc,t^{1/3}b_m$ and $\mathcal{C}_t'=\mathcal{C}_t+\frac12t^{1/3}$.
We scale $(\eta,\xi)\mapsto(t^{1/3}\eta,t^{1/3}\xi)$ and $(u,v)\mapsto(t^{-1/3}u,t^{-1/3}v)$ (in the Fredholm determinant) so that $\fK_0(u,v)$ is now given as $\int_{\mathcal{C}_1} d\eta\int_{\mathcal{C}_1'} d\xi\,\frac{\pi}{\sin(\pi(\xi-\eta))}\frac{e^{t\xi^3/3-(u+r)xi}}{e^{t\eta^3/3-(v+r)\eta}}\prod_{k=1}^m\!\frac{\Gamma(\eta-b_k)}{\Gamma(\xi-b_k)}$.
Since $\kh_{b+x/t}(t,x)+x^2/t$ is stationary in $x$ (see \cite[Rem. 1.14]{borodinCorwinFerrari}), we may write $G_{b}(t,x,r)=\det(\fI-\fK)_{L^2[0,\infty)}$ with
\begin{equation}
\fK(u,v)=\itwopii{2}\int_{\widetilde{\mathcal{C}}_1} d\xi\int_{\widetilde{\mathcal{C}}_1'} d\eta\,\frac{\pi}{\sin(\pi(\xi-\eta))}\frac{e^{t\xi^3/3-(u+r+x^2/t)\xi}}{e^{t\eta^3/3-(v+r+x^2/t)\eta}}\prod_{k=1}^m\frac{\Gamma(\eta-b_k-x/t)}{\Gamma(\xi-b_k-x/t)}
\end{equation}
(here the contour $\widetilde{\mathcal{C}}_1$ has to cross the real axis to the right of $b_i+x/t$ for all $i$).
Conjugating the kernel by $e^{ux/t}$ as in the previous case (i.e. multiplying the kernel by $e^{(u-v)x/t}$) and changing variables $\eta\longmapsto\eta+x/t$, $\xi\longmapsto\xi+x/t$, we get
\begin{equation}\label{eq:gumbel-spiked-fred}
G_{b}(t,x,r)=\det(\fI-\wt\fK)_{L^2[0,\infty)}
\end{equation}
with
\[\wt\fK(u,v)=\itwopii{2}\int_{\mathcal{C}_1} d\xi\int_{\mathcal{C}_1'} d\eta\,\frac{\pi}{\sin(\pi(\xi-\eta))}\frac{e^{t\xi^3/3+x\xi^2-(u+r)\xi}}{e^{t\eta^3/3+x\eta^2-(v+r)\eta}}\prod_{k=1}^m\frac{\Gamma(\eta-b_k)}{\Gamma(\xi-b_k)}\]
As in the previous case, $\tilde\fK$ satisfies the necessary differential relations, so it follows again that:

\begin{thm}\label{thm:spiked}
 $\phi_b\coloneqq \p_r^2\log G_b$ satisfies the KP-II equation \eqref{eq:KP-II}
\end{thm}

\subsection*{KPZ equation with two-sided Brownian initial data}

In \cite{ledoussal-kp} (which appeared about a month after the first version of the present article), P. Le Doussal suggested that a certain modified generating function for the solution $\kh_{w_\pm}$ of \eqref{KPZ} with the same scaling as above and initial data of the form $\kh_{w_\pm}(0,x)=\fB(x)+w_-x\uno{x<0}+w_+x\uno{x\geq0}$ with $\fB$ a double-sided Brownian motion with $\fB(0)=0$ and $w_->w_+$ should also satisfy KP.
This is actually true, as we explain next.
Define
\[\wt G_{w_\pm}(t,x,r)=\Gamma(w_--w_+)^{-1}\ee\!\left[2\tts e^{\frac12(w_--w_+)(\kh_{w_\pm}(t,x)+\frac{t}{12}-r)}K_{w_+-w_-}(2e^{\frac12(\kh_{w_\pm}(t,x)+\frac{t}{12}-r)})\right]\]
where $K_\nu$ is the modified Bessel function of order $\nu$.
This modified generating function $\wt G_{w_\pm}$ can alternatively be expressed as the analog of \eqref{eq:gumbel} where the KPZ height function $\kh_{w_\pm}$ is replaced by a randomly shifted height function $\kh_{w_\pm}(t,x)+\Upsilon$,
\[\wt G_{w_\pm}(t,x,r)=\ee\!\left[\exp\!\left\{-e^{\kh_{w_\pm}(t,x)+\Upsilon+\frac{t}{12}-r}\right\}\right]\]
with $\Upsilon$ an independent log-gamma random variable with parameter $w_--w_+$, i.e. $e^{-\Upsilon}$ has density $\Gamma(w_--w_+)^{-1}x^{w_--w_+-1}e^{-x}$.
Explicit formulas for the distribution of this shifted height function were obtained in \cite{imamSasamStationary,imamSasamStat2} using the non-rigorous replica method.
A similar formula, which is the one we will use below, was obtained rigorously later on \cite{BCFV}; the equality between the two expressions above for $\wt G_{w_\pm}$ is essentially Cor. 2.6 in the latter paper, see also Rem. 2.10 there.

From \cite[Thm. 2.9]{BCFV} (and footnote \ref{ft:sc} again) we have $\wt G_{w_\pm}(t,x,r)=\det(\fI-\fK)_{L^2[0,\infty)}$ with
\begin{multline}
\fK(u,v)=\itwopii{2}\int_{\mathcal{C}_t} d\eta\int_{\tilde{\mathcal{C}}_t} d\xi\,\frac{t^{-1/3}\pi}{\sin(t^{-1/3}\pi(\xi-\eta))}\frac{e^{\xi^3/3-(u+t^{-1/3}(r+\frac{x^2}{t}))\xi}}{e^{\eta^3/3-(v+t^{-1/3}(r+\frac{x^2}{t}))\eta}}\\
\times\frac{\Gamma(w_-+\frac{x}{t}-t^{-1/3}\xi)}{\Gamma(t^{-1/3}\xi-w_+-\frac{x}{t})}\frac{\Gamma(t^{-1/3}\eta-w_+-\frac{x}{t})}{\Gamma(w_-+\frac{x}{t}-t^{-1/3}\eta)}
\end{multline}
where $\mathcal{C}_t$ goes from $-\frac14t^{1/3}-\I\infty$ to $-\frac14t^{1/3}+\I\infty$ crossing the real axis between $t^{1/3}w_+$ and $t^{1/3}w_-$ and $\tilde{\mathcal{C}}_t$ goes from $\frac14t^{1/3}-\I\infty$ to $\frac14t^{1/3}+\I\infty$ staying to the right of $\mathcal{C}_t$ and also crossing the real axis between $t^{1/3}w_+$ and $t^{1/3}w_-$.
Changing variables as in the last example leads to
\begin{equation}\label{eq:2BMgenfn}
\wt G_{w_\pm}(t,x,r)=\det(\fI-\wt\fK)_{L^2[0,\infty)}
\end{equation}
with
\[\wt\fK(u,v)=\itwopii{2}\int_{\mathcal{C}_1} d\eta\int_{\tilde{\mathcal{C}}_1} d\xi\,\frac{\pi}{\sin(\pi(\xi-\eta))}\frac{e^{t\xi^3/3+x\xi^2-(u+r)\xi}}{e^{t\eta^3/3+x\eta^2-(v+r)\eta}}\frac{\Gamma(w_--\xi)}{\Gamma(\xi-w_+)}\frac{\Gamma(\eta-w_+)}{\Gamma(w_--\eta)}.\]
As in the previous examples, $\wt\fK$ satisfies the hypotheses of Thm. \ref{thm3}, so:

\begin{thm}
 $\wt\phi_{w_\pm}(t,x,r)\coloneqq\p_r^2\log\wt G_{w_\pm}(t,x,r)$ satisfies the KP-II equation \eqref{eq:KP-II}.
 \end{thm}

\vs
Our next two examples are related with the setting of Thm. \ref{thm:spiked}. 
Here we abandon the setting of the KPZ equation and go back to random matrix distributions and the KPZ fixed point.

\subsection*{The BBP distribution for spiked random matrices}

Consider $G_b$ as in the case of spiked initial data.
It is known (see \cite[Cor. 1.15]{borodinCorwinFerrari} for the case $t=1,x=0$, the general case follows in the same way or by scaling and shift invariance) that
\[\lim_{\ep\to0}G_{\ep^{1/2}b}(\ep^{-3/2}t,\ep^{-1}x,\ep^{-1/2}r)=F_{\uptext{BBP},t^{1/3}b-t^{-2/3}x}(t^{-1/3}r+t^{-4/3}x^2)\]
where $F_{\uptext{BBP},b}$ is the Baik-Ben Arous-P\'ech\'e (BBP) distribution arising from spiked (unitarily invariant) random matrices \cite{BBP}.
On the other hand, by scaling invariance of KP-II (see Rem. \ref{rem:m}), for each fixed $\ep>0$, $\phi^\ep_b(t,x,r)\coloneqq\p_r^2\log G_{\ep^{1/2}b}(\ep^{-3/2}t,\ep^{-1}x,\ep^{-1/2}r)$ satisfies KP-II as well.
As a consequence, one expects that:

\begin{thm}\label{thm:bbp}
Define
\[\tilde F_{\uptext{BBP},b}(t,x,r)=F_{\uptext{BBP},t^{1/3}b-t^{-2/3}x}(t^{-1/3}r+t^{-4/3}x^2).\]
Then $\phi_{\uptext{BBP},b}\coloneqq\p_r^2\log\tilde F_{\uptext{BBP},b}$ satisfies the KP-II equation \eqref{eq:KP-II}.
\end{thm}

This is indeed the case, as can be checked in a similar way as above using Thm. \ref{thm3} and the explicit Fredholm determinant formula for $F_{\uptext{BBP},b}$ (see for instance \cite[Eqn. (1.4)]{borodinCorwinFerrari}).

\subsection*{KPZ fixed point with half-Brownian initial data}
Consider the KPZ fixed point $\fh$ with half-Brownian initial data $\fh_0(x)=\fB(-x)$ for $x\leq0$ and $\fh_0(x)=-\infty$ for $x>0$, where $\fB$ is a Brownian motion started at $0$ with diffusivity $2$.
In the setting of \cite{fixedpt}, this corresponds to starting TASEP with a product measure with density $\frac12$ on the negative integers and no particles on the positive integers, and known results \cite{corwinFerrariPeche} in this case give $\pp(\fh(1,x)\leq r)=F_{\uptext{BBP},x}(r+x^2)$.
By scaling invariance of the KPZ fixed point, this gives in this case for $F_\uptext{half-BM}(t,x,r)=\pp_{\fh_0}(\fh(t,x)\leq r)$
\[F_\uptext{half-BM}(t,r,x)=F_{\uptext{BBP},t^{-2/3}x}(t^{-1/3}r+t^{-4/3}x^2).\]
Comparing with Thm. \ref{thm:bbp} we deduce that:

\begin{thm}
 $\phi_\uptext{half-BM}=\p_r^2\log F_\uptext{half-BM}$ satisfies the KP-II equation \eqref{eq:KP-II}.
\end{thm}

Similar statements can be written for more general spiked initial conditions in this setting as well as for multipoint distributions (now in terms of the matrix KP equation as in Thm. \ref{thm:2}).

\settocdepth{subsection}
\hypersetup{bookmarksdepth=2}

\section{Derivation of KP}\label{sec:proof}

This section contains the proof of Thm. \ref{thm3}.
After we performed the complicated computation, we discovered that a very similar argument was actually known in the one dimensional case \cite{poppeIP}. 
It is shown there that the Fredholm determinant of a kernel satisfying suitable differential relations solves the Hirota equations.
The differential relations turn out to be equivalent to the way the kernel
depends on $t$, $x$ and $r$ in Thm. \ref{thm3}.
It seems to actually go back to \cite{zaharovShabat2,zaharovShabat1} though it is not explicit there, and rediscovered in the literature multiple times.

\subsection{The logarithmic derivative}\label{sec:thirdsec}

Let
\[\Phi(t,r,x)=\p_r\log\det(\fI-\fK)_{L^2(\rr^+)\oplus\dotsm\oplus L^2(\rr^+)}\]
where $\fK$ is a matrix kernel satisfying the assumptions of Thm. \ref{thm3}.
From now on we omit the subscript on the Fredholm determinant and traces.

Given an operator $\fA$ acting on the $n$-fold direct sum of $L^2(\rr^+)$ with matrix kernel $(\fA_{ab}(u,v))_{a,b=1}^n$ we will write $\D_1\fA$ and $\D_2\fA$ for the operators with kernels given by
\[\D_1\fA(u,v)=\big(\p_u\fK_{ab}(u,v)\big)_{a,b=1,\dotsc,n}\qand\D_2\fA(u,v)=\big(\p_v\fK_{ab}(u,v)\big)_{a,b=1,\dotsc,n}.\]
Note that this is just a notational device; $\D_i$ is not meant to denote an operator.
By (1) of Thm. \ref{thm3},
\begin{align}
\Phi(t,r,x)&=\p_{r}\tsm\log(\det(\fI-\fK))=-\tr((\fI-\fK)^{-1}\p_{r}\fK)\\
&=-\tr((\fI-\fK)^{-1}(\D_1+\D_2)\fK)=-\tr((\fI-\fK)^{-1}\D_1\fK+\D_2((\fI-\fK)^{-1}\fK))\\
&=-\tr(\D_1\fK(\fI-\fK)^{-1}+\D_2((\fI-\fK)^{-1}\fK))=-\tr((\D_1+\D_2)((\fI-\fK)^{-1}\fK))\\
&=-\sum_a\int_0^\infty d\xi\,\p_\xi\big((\fI-\fK)^{-1}\fK\big)_{aa}(\xi,\xi)=\sum_a\big((\fI-\fK)^{-1}\fK\big)_{aa}(0,0),
\end{align}
where we used the cyclicity of the trace.
Introducing the notation
\[[\fA]=\big(\fA_{a,b}(0,0)\big)_{a,b=1,\dotsc,n},\]
this tells us that $\Phi$ can be expressed as an $n$-dimensional trace,
\[\Phi(t,r,x)=\tr[\fR\fK]\qquad\uptext{with}\quad\fR=(\fI-\fK)^{-1}.\]
Note here that $\fI-\fK$ is invertible because the determinant is non-zero.  This proves \eqref{traceQ} (after identifying $\p_r$ with $\cD_r$ as in Sec. \ref{sec:main}) with, in the present notation,
\begin{equation}\label{eq:defQ}
Q=[\fR\fK].
\end{equation}

\subsection{Formulas for the partial derivatives}

The goal is now is to use (1), (2), (3) of Thm. \ref{thm3} to show that algebraically $Q$ satisfies the matrix KP equation \eqref{eq:matKP}.

Write $\fK'=\p_r\fK$.
Using the general formula $\p_a(\fI-\fA(a))^{-1}=(\fI-\fA(a))^{-1}\p_a\fA(a)(\fI-\fA(a))^{-1}$ for an operator $\fA(a)$ depending smoothly on a parameter $a$ together with the identity
\[\fK\fR=\fR\fK=\fR-\fI,\]
which we will use repeatedly, we have
\begin{align}
\p_rQ&=[\p_r(\fR\fK)]=[\fR\fK'\fR\fK+\fR\fK']=[\fR\fK'\fR],\\
\shortintertext{and similarly}
\p_r^2Q&=2[\fR\fK'\fR\fK'\fR]+[\fR\fK''\fR],\\
\p_r^3Q&=6[\fR\fK'\fR\fK'\fR\fK'\fR]+3[\fR\fK''\fR\fK'\fR]+3[\fR\fK'\fR\fK''\fR]+[\fR\fK'''\fR],\\
\p_tQ&=[\fR\p_t\fK\fR].
\end{align}

Next we want to compute $(\p_rQ)^2$.
Note that, in general,
\begin{align}
([\fA][\fBb])_{a,b}&=-\sum_{c}\int_0^{\infty}d\eta\,\p_\eta(\fA_{ac}(0,\eta)\fBb_{cb}(\eta,0))\\
&=-\sum_{c}\int_0^{\infty}d\eta\,\big(\p_\eta\fA_{ac}(0,\eta)\fBb_{cb}(\eta,0)+\fA_{ac}(0,\eta)\p_\eta\fBb_{cb}(\eta,0)\big)\\
&=-(\D_2\fA\fBb)_{ab}(0,0)-(\fA\D_1\fBb)_{ab}(0,0)
\end{align}
so that the following integration by parts formula holds:
\begin{equation}
[\fA][\fBb]=-[\fA\D_1\fBb+\D_2\fA\fBb].\label{eq:multi-ibp}
\end{equation}
We will use this in the formula
\[(\p_rQ)^2=([\fR\fK'\fR\fK]+[\fR\fK'])^2 
=[\fR\fK'\fR\fK]^2+[\fR\fK'\fR\fK][\fR\fK']+[\fR\fK'][\fR\fK'\fR\fK]+[\fR\fK']^2.\]
The first term equals (using $\D_2(\fK_a\fK_b)=\fK_a\D_2\fK_b$)
\[-[\fR\fK'\fR(\D_2\fK\fR+\fK\D_1\fR)\fK'\fR\fK]=-[\fR\fK'\fR\fK'\fR\fK'\fR\fK]+[\fR\fK'\fR(\D_1\fK\fR-\fK\D_1\fR)\fK'\fR\fK].\]
Similarly, the fourth term equals
\[-[\fR\fK''\fR\fK']+[\fR(\D_1\fK'\fR-\fK'\D_1\fR)\fK']\]
and the two middle ones equal
\[-[\fR\fK'\fR\fK'\fR\fK']-[\fR\fK''\fR\fK'\fR\fK]+[\fR\fK'\fR(\D_1\fK\fR-\fK\D_1\fR)\fK']+[\fR(\D_1\fK'\fR-\fK'\D_1\fR)\fK'\fR\fK].\]
Using this together with (1) and (2) of Thm. \ref{thm3}, which in the current notation read $\p_r\fK=(\D_1+\D_2)\fK$ and $\p_t\fK=-\frac13(\D_1^3+\D_2^3)\fK$, we get
\begin{align}
12\p_t&Q(t,r,x)+\p_r^3Q(t,r,x)+6(\p_rQ(t,r,x))^2\\
&=-4[\fR(\D_1^3+\D_2^3)\fK\fR]+6[\fR\fK'\fR\fK'\fR\fK'\fR]+3[\fR\fK''\fR\fK'\fR]\label{123}\tag*{$-$(1)$+$(2)$+$(3)}\\
&\hspace{0.5in}+3[\fR\fK'\fR\fK''\fR]+[\fR(\D_1^3+\D_2^3+3\D_1^2\D_2+3\D_1\D_2^2)\fK\fR]\label{4}\tag*{$+$(4)$+$(5)}\\
&\hspace{0.35in}-6[\fR\fK''\fR\fK']+6[\fR(\D_1\fK'\fR-\fK'\D_1\fR)\fK']\label{56}\tag*{$-$(6)$+$(7)}\\
&\hspace{0.35in}-6[\fR\fK'\fR\fK'\fR\fK']-6[\fR\fK''\fR\fK'\fR\fK]\tag*{$-$(8)$-$(9)}\\
&\hspace{0.5in}\label{1011}\tag*{$+$(10)$+$(11)}+6[\fR\fK'\fR(\D_1\fK\fR-\fK\D_1\fR)\fK']+6[\fR(\D_1\fK'\fR-\fK'\D_1\fR)\fK'\fR\fK]\\
&\hspace{0.35in}-6[\fR\fK'\fR\fK'\fR\fK'\fR\fK]+6[\fR\fK'\fR(\D_1\fK\fR-\fK\D_1\fR)\fK'\fR\fK{}].\label{1213}\tag*{$-$(12)$+$(13)}
\end{align}
(12) equals $6[\fR\fK'\fR\fK'\fR\fK'\fR]-6[\fR\fK'\fR\fK'\fR\fK']$, so
\begin{align}
S_{\uptext{I}}&\coloneqq-\uptext{(1)}+\uptext{(2)}+\uptext{(5)}-\uptext{(8)}-\uptext{(12)}=-3[\fR(\D_1^3+\D_2^3-\D_1^2\D_2-\D_1\D_2^2)\fK\fR]\\
&=-3[\fR(\D_1-\D_2)^2(\D_1+\D_2)\fK\fR]=-3[\fR(\D_1-\D_2)^2\fK'\fR].
\end{align}
Similarly
\begin{equation}
S_{\uptext{II}}\coloneqq+\uptext{(3)}+\uptext{(4)}-\uptext{(6)}-\uptext{(9)}=-3[\fR\fK''\fR\fK'\fR]+3[\fR\fK'\fR\fK''\fR]
\end{equation}
and
\begin{equation}
S_{\uptext{III}}\coloneqq+\uptext{(7)}+\uptext{(10)}+\uptext{(11)}+\uptext{(13)}
=+6[\fR(\D_1\fK'\fR-\fK'\D_1\fR)\fK'\fR]+6[\fR\fK'\fR(\D_1\fK\fR-\fK\D_1\fR)\fK'\fR].
\end{equation}
So
\begin{equation}
\begin{aligned}
&4\p_tQ(t,r,x)+\tfrac13\p_r^3Q(t,r,x)+2(\p_rQ(t,r,x))^2=\tfrac13(S_{\uptext{I}}+S_{\uptext{II}}+S_{\uptext{III}})\\
&\hspace{0.5in}=-[\fR\fK''\fR\fK'\fR]+[\fR\fK'\fR\fK''\fR]-[\fR(\D_1^3+\D_2^3-\D_1^2\D_2-\D_1\D_2^2)\fK\fR]\\
&\hspace{1in}+2[\fR(\D_1\fK'\fR-\fK'\D_1\fR)\fK'\fR]+2[\fR\fK'\fR(\D_1\fK\fR-\fK\D_1\fR)\fK'\fR]\\
&\hspace{0.5in}=-[\fR(\D_2-\D_1)\fK'\fR\fK'\fR]+[\fR\fK'\fR\fK''\fR]-[\fR(\D_1^3+\D_2^3-\D_1^2\D_2-\D_1\D_2^2)\fK\fR]\\
&\hspace{1in}-2[\fR\fK'\fR\D_1\fR\fK'\fR]+2[\fR\fK'\fR\D_1\fK\fR\fK'\fR].
\end{aligned}\label{eq:firsteq}
\end{equation}
Using now $-[\fR\fK'\fR\D_1\fR\fK'\fR]+[\fR\fK'\fR\D_1\fK\fR\fK'\fR]=-[\fR\fK'\fR\D_1(\fI-\fK)\fR\fK'\fR]$, which equals $-[\fR\fK'\fR\D_1\fK'\fR]$, yields
\begin{multline}\label{eq:ikdvterm}
4\p_tQ(t,r,x)+\tfrac13\p_r^3Q(t,r,x)+2(\p_rQ(t,r,x))^2\\
=-[\fR(\D_2-\D_1)\fK'\fR\fK'\fR]-[\fR\fK'\fR(\D_1-\D_2)\fK'\fR]-[\fR(\D_1-\D_2)^2\fK'\fR].
\end{multline}

\vspace{2pt}
\begin{rem}
At this stage we can already see that the one point distribution in the flat case $\fh_0\equiv0$ satisfies the (integrated) KdV equation.
In fact, the arguments in \cite[Sec. 4.4]{fixedpt} lead in this case (in the language of the Brownian scattering operator \eqref{sg}) to $\fK^{\hypo(\fh)}_t=\fI-\fU_t(\fI-\varrho)\bar\P_0(\fI-\varrho)\fU_{t}^{-1}=\fU_t\varrho\tts\fU_{t}^{-1}$ with $\varrho f(x)=f(-x)$, which yields $\fK(u,v)=(2t)^{-1/3}\Ai((2t)^{-1/3}(u+v+2r))$.
Hence $\fK$ is a Hankel kernel, and thus the right hand side in \eqref{eq:ikdvterm} vanishes.
\end{rem}

Next we add the derivatives in the $x_i$ variables.
As for $\p_r^2Q$, we have
\begin{equation}
\p_x^2Q=2[\fR\p_x\fK\fR\p_x\fK\fR]+[\fR\p_x^2\fK\fR].\label{eq:eq:Dx2Q}
\end{equation}
On the other hand, if we apply $\p_r$ to \eqref{eq:ikdvterm} and use (3) of Thm. \ref{thm3} to write $(\D_2-\D_1)\fK'=\p_x\fK$ we get
\begin{equation}\label{eq:kdvterm}
\begin{split}
&-[\fR\fK'\fR\p_x\fK\fR\fK'\fR]-[\fR\p_x\fK'\fR\fK'\fR]-2[\fR\p_x\fK\fR\fK'\fR\fK'\fR]-[\fR\p_x\fK\fR\fK''\fR]\\
&\hspace{0.4in}+2[\fR\fK'\fR\fK'\fR\p_x\fK\fR]+[\fR\fK''\fR\p_x\fK\fR]+[\fR\fK'\fR\p_x\fK'\fR]+[\fR\fK'\fR\p_x\fK\fR\fK'\fR]\\
&\hspace{0.4in}+[\fR\fK'\fR(\D_1-\D_2)\p_x\fK\fR]+[\fR(\D_1-\D_2)\p_x\fK'\fR]+[\fR(\D_1-\D_2)\p_x\fK\fR\fK'\fR].
\end{split}
\end{equation}
Note that the first and eighth terms cancel.
We want to add $\p_x^2Q(t,r,x)$.
Since $(\D_1-\D_2)\p_x\fK'=\p_x(\D_1^2-\D_2^2)\fK=-\p_x^2\fK$, the next-to-last term in the last expression cancels the second bracket on the right hand side of \eqref{eq:eq:Dx2Q}.
Using additionally $\p_x\fK'+(\D_1-\D_2)\p_x\fK=2\D_1\p_x\fK$ and $-\p_x\fK'+(\D_1-\D_2)\p_x\fK=-2\D_2\p_x\fK$  and writing
\begin{equation}
q=\p_rQ,
\end{equation}
we deduce that
\begin{equation}\label{eq:full1-multi}
\begin{split}
&4\p_tq+\tfrac13\p_r^3q+2(q\tts\p_rq+\p_rq\ts q)+\p_x^2Q(t,r,x)\\
&\hspace{0.35in}=2\ts\Big(\tsm\!-[\fR\D_2\p_x\fK\fR\fK'\fR]-[\fR\p_x\fK\fR\fK'\fR\fK'\fR]-\tfrac12[\fR\p_x\fK\fR\fK''\fR]\\
&\hspace{1.4in}+[\fR\fK'\fR\fK'\fR\p_x\fK\fR]+\tfrac12[\fR\fK''\fR\p_x\fK\fR]+[\fR\fK'\fR\D_1\p_x\fK\fR]\\
&\hspace{1.4in}+[\fR\p_x\fK\fR\p_x\fK\fR]\Big).
\end{split}
\end{equation}
We claim that the right hand side equals two times
\begin{equation}\label{eq:full1'}
\begin{split}
&-\big( [\fR\p_x\fK\fR\D_1\fK\fR\fK'\fR] + [\fR\p_x\D_2\fK\fR\fK'\fR]\big)
+ \big( [\fR\fK'\fR\D_2\fK\fR\p_x\fK\fR] + [\fR\fK'\fR\D_1\p_x\fK\fR]\big)\\
&-\big( [\fR\p_x\fK\fR\D_2\fK\fR\fK'\fR] + [\fR\p_x\fK\fR \D_1\fK'\fR]\big)
 + \big( [\fR\fK'\fR\D_1\fK\fR\p_x\fK\fR] + [\fR\D_2\fK'\fR\p_x\fK\fR]\big). 
\end{split}
\end{equation}
To see this, express the right hand side of \eqref{eq:full1-multi} as $2(r_1+r_2+\dotsc+r_7)$, express \eqref{eq:full1'} as $q_1+q_2+\dotsc+q_8$, and note first that $r_1=q_2$, $r_6=q_4$, $r_2=q_1+q_5$ and $r_4=q_3+q_7$.
On the other hand we have $r_3=-\frac12[\fR\p_x\fK\fR(\D_1+\D_2)\fK'\fR]=q_6+\frac12[\fR\p_x\fK\fR(\D_1-\D_2)\fK'\fR]=q_6-\frac12r_7$ and similarly $r_5=\frac12[\fR(\D_1+\D_2)\fK'\fR\p_x\fK\fR]=q_8+\frac12[\fR(\D_1-\D_2)\fK'\fR\p_x\fK\fR]=q_8-\frac12r_7$.
This gives $r_3+r_5+r_7=q_6+q_8$, and finishes proving the claim.

Integrating by parts (i.e. using \eqref{eq:multi-ibp}) within each parenthesis in \eqref{eq:full1'} we get
\begin{equation}
\begin{split}
2\p_tq+\tfrac16\p_r^3q+(q\tts\p_rq&+\p_rq\ts q)+\tfrac12\p_x^2Q(t,r,x)\\
=&-[\fR\p_x\fK(\fR\D_1\fK\fR-\D_1\fR) \fK'\fR] +  [\fR\p_x\fK][\fR\fK'\fR]\\
& + [\fR\fK'(\fR\D_2\fK\fR-\D_2\fR)\p_x\fK\fR] - [\fR\fK'\fR][\p_x\fK\fR]\\
&- [\fR\p_x\fK(\fR\D_2\fK\fR-\D_2\fR)\fK'\fR] +  [\fR\p_x\fK\fR][\fK'\fR] \\
& +   [\fR\fK'(\fR\D_1\fK\fR-\D_1\fR)\p_x\fK\fR]  - [\fR\fK'][\fR\p_x\fK\fR].
\end{split}
\end{equation}
Write this as $s_1+\dotsm+s_8$.
Notice that $s_1+s_5$ yields a term involving
\[( \fR\D_1\fK \fR-\D_1\fR)+ (\fR\D_2\fK\fR-\D_2\fR) = \fR(\fK\D_1\fK+\D_2\fK\fK)\fR-(\D_1+\D_2)\fI,\]
where we have used $\fR\fK=\fK\fR=\fR-\fI$ again, and thus integrating by parts one more time we get 
\begin{align}
s_1+s_5&=-[\fR\p_x\fK(\fR(\fK\D_1\fK+\D_2\fK\fK)\fR-(\D_1+\D_2)\fI)\fK'\fR]\\
&=[\fR\p_x\fK\fR\fK][\fK\fR\fK'\fR]+[\fR\p_x\fK((\D_1+\D_2)\fI)\fK'\fR]\\
&=([\fR\p_x\fK\fR]-[\fR\p_x\fK])([\fR\fK'\fR]-[\fK'\fR])+[\fR\p_x\fK((\D_1+\D_2)\fI)\fK'\fR]\\
&=[\fR\p_x\fK\fR][\fR\fK'\fR]-s_2-s_6+[\fR\p_x\fK][\fK'\fR]+[\fR\p_x\fK((\D_1+\D_2)\fI)\fK'\fR].
\end{align}
In a similar fashion we get
\[s_3+s_7=-[\fR\fK'\fR][\fR\p_x\fK\fR]-s_4-s_8-[\fR\fK'][\p_x\fK\fR]-[\fR\fK'((\D_1+\D_2)\fI)\p_x\fK\fR].\]
Therefore
\begin{multline}
2\p_tq+\tfrac16\p_r^3q+(q\tts\p_rq+\p_rq\ts q)+\tfrac12\p_x^2Q(t,r,x)
=[\fR\p_x\fK][\fK'\fR]+[\fR\p_x\fK((\D_1+\D_2)\fI)\fK'\fR]\\
-[\fR\fK'][\p_x\fK\fR]-[\fR\fK'((\D_1+\D_2)\fI)\p_x\fK\fR]+[\fR\p_x\fK\fR][\fR\fK'\fR]-[\fR\fK'\fR][\fR\p_x\fK\fR].
\end{multline}
In order to complete the proof we note that if $\fA$ and $\fBb$ are nice kernels then integrating by parts we get
\begin{equation}\label{eq:trick2}
[\fA((\D_1+\D_2)\fI)\fBb]=[\fA\D_1\fBb]+[\fA\D_2\fI\tts\fBb]=-[\fA][\fBb],
\end{equation}
which immediately yields 
\begin{equation}
2\p_tq+\tfrac16\p_r^3q+(q\tts\p_rq+\p_rq\ts q)+\tfrac12\p_x^2Q(t,r,x)=[\fR\p_x\fK\fR][\fR\fK'\fR]-[\fR\fK'\fR][\fR\p_x\fK\fR].
\end{equation}
The right hand side equals $\p_xQ\ts q-q\tts\p_xQ$, so the equation becomes \eqref{eq:matKP} as needed.

\subsection{Proof of Thm. \ref{thm3}}

For $r$ inside the interval where the trace norm of $\fK$ is strictly less than $1$ everything in the previous section is well defined, because that condition ensures that $(\fI-\fK)^{-1}-\fI$ is trace class, and therefore the series for the brackets are convergent.
So the algebraic computations hold pointwise for such $r$.    

Since the kernel $\fK(u,v)$ is real analytic in $t>0$, $x$ and $r$,
the Fredholm determinant $\det(\fI-\fK(t,x,r))$ is as well, since it is given by its  Fredholm series, each of whose terms is real analytic, and which converges uniformly since $\fK$ is trace class uniformly in compact sets of $t,x,r$.
Since the determinant never vanishes, $\p^2_r \log \det(\fI-\fK)$ is also real analytic, as are all the terms in the KP equation (either the scalar or matrix version).
Therefore the left hand side of the KP equation \eqref{eq:matKP} (or \eqref{eq:KP-II})
 is a real analytic function.
We have proved this function vanishes for $r$ in an open interval, and therefore it vanishes everywhere.

\appendix

\section{Multipoint initial data}
\label{mid}

\subsection{\texorpdfstring{$t\to0$}{t to 0} limit of the Brownian scattering operator}

Let the initial data for the KPZ fixed point be a finite collection of narrow wedges $\mathfrak{d}^{\vec b}_{\vec a}$ as in Ex. \ref{fcnw}.
Fix $x_1<\dotsc<x_m$.
Our goal is to compute the limit
\begin{equation}\label{eq:icker}
\lim_{t\to0}e^{-x_i\p^2}\fK^{\hypo(\mathfrak{d}_{\vec a}^{\vec b})}_te^{x_j\p^2}
\end{equation}
in operator norm in $L^2([r,\infty))$ for any fixed $r$, with $\fK^{\hypo(\fh)}_t$ defined in \eqref{sg}.
The convergence could be upgraded to trace norm in most cases, but it does not hold, for example, in the important case $m=1$, $x_1=a_1$, where it is easy to see that the limit itself is not trace class (see \eqref{eq:xiaxjcase}).
Nor should it be; one should not expect convergence to the initial data in such a strong sense, what we are aiming for here is just to understand the $t\to0$ behavior of the matrix kernel.

Recall the operators $\fT_{t,x}$ introduced in \eqref{eq:fTdef}, which satisfy the differential relations \eqref{eq:partial-fT}.
They can be expressed as $\fT_{t,x}=e^{-\frac13t\p^3+x\p^2}=e^{x\p^2}\fU_t$, which means that 
\[\fK^{\hypo(\fh)}_t=\lim_{\ell\to\infty}\fT_{-t,-\ell}{\bf P}_{-\ell,\ell}^{\uptext{Hit}\,\fh}\fT_{t,-\ell}.\]


We begin by studying the single narrow wedge case, $k=1$, writing $a=a_1$, $b=b_1$.
In this case we have, for $\ell>a$, ${\bf P}_{-\ell,\ell}^{\uptext{Hit}\,\mathfrak{d}_a^b}=e^{(a+\ell)\p^2}\bar\P_be^{(\ell-a)\p^2}$ by definition of the left hand side and where $\bar\P_b$ was defined in \eqref{eq:Pdef}.
Then
\begin{equation}\label{eq:nwKt}
\fK^{\hypo(\mathfrak{d}_{a}^{b})}_t
=\lim_{\ell\to\infty}\fT_{-t,-\ell}e^{(a+\ell)\p^2}\bar\P_be^{(\ell-a)\p^2}\fT_{t,-\ell}=\fT_{-t,a}\bar\P_b\fT_{t,-a}.
\end{equation}
Consider first the case $x_i\leq a\leq x_j$.
Since $\fT_{\pm t,x}\longrightarrow e^{x\p^2}$ as $t\to0$ for all $x\geq0$ in operator norm, by \eqref{eq:nwKt} we get
\begin{equation}\label{eq:xiaxjcase}
\lim_{t\to0}e^{-x_i\p^2}\fK^{\hypo(\mathfrak{d}_{a}^{b})}_te^{x_j\p^2}=\lim_{t\to0}\fT_{-t,a-x_i}\bar\P_b\fT_{t,x_j-a}
=e^{(a-x_i)\p^2}\bar\P_be^{(x_j-a)\p^2}={\bf P}^{\uptext{Hit }\mathfrak{d}_{a}^{b}}_{x_i,x_j}
\end{equation}
in operator norm and in all of $L^2(\rr)$.

Next consider the case  $x_i\leq x_j$ and $a\notin[x_i,x_j]$.
We will show that in this case our operator goes to $0$ as $t\to0$.
We will assume for simplicity that $a>x_j$, the case $a<x_i$ works in the same way.
We have $\P_re^{-x_i\p^2}\fK^{\hypo(\mathfrak{d}_{a}^{b})}_te^{x_j\p^2}\P_r=(\P_r\fT_{-t,a-x_i}\bar\P_b)(\bar\P_b\fT_{t,x_j-a}\P_r)$, and the first factor goes to $\P_r e^{(a-x_i)\p^2}\bar\P_b$ as $t\to0$ in operator norm, so it is enough to show that $\bar\P_b\fT_{t,x_j-a}\P_r$ goes to $0$.
We estimate its Hilbert-Schmidt norm,
\begin{equation}\label{eq:HS}
\|\bar\P_b\fT_{t,x_j-a}\P_r\|_2^2=\int_r^\infty dv\int_{-\infty}^b d\eta\,t^{-2/3}e^{4\bar x_j^3/3t^2-2(\eta-v)\bar x_j/t}\tsm\ttsm\Ai(t^{-1/3}(v-\eta)+t^{-4/3}\bar x_j^2)^2,
\end{equation}
where $\bar x_j=x_j-a$.
Split the $\eta$ integral according to whether $\eta\leq v\wedge b$ or $v\wedge b<\eta\leq b$.
On the first piece we may use the classical bound on the Airy function $|\tsm\ttsm\Ai(s)|\leq C\tts e^{-\frac23(s\vee0)^{3/2}}$ to see that the integral is bounded by $C\tts t^{-2/3}e^{\frac43\bar x_j^3/t^2}\int_r^\infty dv\int_{-\infty}^{v\wedge b}d\eta\,e^{-\frac43(v-\eta)^{3/2}/t^{1/2}-\frac43|\bar x_j|^3/t^2-2(\eta-v)\bar x_j/t}$.
The exponent in the $\eta$ integral is maximized at the edge of the integration, $\eta=v\wedge b$, so applying Laplace's method we deduce that the same integral is bounded by
$C\tts t^{-c}\tts e^{-\frac83|\bar x_j|^3/t^2}\tsm\left[|r-b|+\int_b^\infty dv\,e^{-\frac43(v-b)^{3/2}/t^{1/2}-2(b-v)\bar x_j/t}\right]$, for some $c,C>0$, which clearly goes to $0$ as $t\to0$.
On the second piece the integration region is bounded, so we get directly (since $\bar x_j<0$) that the integral goes to $0$.
This shows that the left hand side of \eqref{eq:HS} goes to $0$ as $t\to0$, as desired.

The last possibility is that $x_i>x_j$.
In this case the operator goes to $0$ again as $t\to0$.
Now one has to estimate the Hilbert-Schmidt norm of the whole operator $e^{-x_i\p^2}\fK^{\hypo(\mathfrak{d}_a^b)}_t e^{x_j\p^2}=\fT_{-t,a-x_i}\bar\P_b\fT_{t,x_j-a}$ on $L^2([r,\infty))$; the estimates are a bit more tedious but very similar to the ones we used in the last case, so we skip the details.

The conclusion of all this is that, in the case of narrow wedge initial data $\mathfrak{d}_{a}^{b}$,
\begin{equation}
\begin{split}
\lim_{t\to0}e^{-x_i\p^2}\fK^{\hypo(\mathfrak{d}_a^b)}_t e^{x_j\p^2}&=\lim_{t\to0}\fT_{-t,a-x_i}\bar\P_b\fT_{t,x_j-a}\\
&=e^{(a-x_i)\p^2}\bar\P_b e^{(x_j-a)\p^2}\uno{x_i\leq a\leq x_j}={\bf P}^{\uptext{Hit }\mathfrak{d}_a^b}_{x_i,x_j}\uno{x_i\leq x_j}
\end{split}\label{eq:Knw0}
\end{equation}
in operator norm in $L^2([r,\infty))$.

Now we turn to the general case $\fh=\mathfrak{d}_{\vec a}^{\vec b}$.
For $\ell>|a_1|\vee|a_k|$ we have, by inclusion-exclusion,
\begin{equation}
{\bf P}_{-\ell,\ell}^{\uptext{Hit}\,\mathfrak{d}_{\vec a}^{\vec b}}=\sum_{n=1}^k(-1)^{n+1}\!\!\sum_{1\leq p_1<\dotsm<p_n\leq k}\!\!e^{(a_{p_1}-\ell)\p^2}\bar\P_{b_{p_1}}e^{(a_{p_2}-a_{p_1})\p^2}\bar\P_{b_{p_2}}\dotsm e^{(a_{p_n}-a_{p_{n-1}})\p^2}\bar\P_{b_{p_n}}e^{(\ell-a_{p_n})\p^2},
\end{equation}
so proceeding as in \eqref{eq:nwKt} we get that $e^{-x_i\p^2}\fK^{\hypo(\mathfrak{d}_{\vec a}^{\vec b})}_te^{x_j\p^2}$ equals
\begin{equation}
\sum_{n=1}^k(-1)^{n+1}\sum_{1\leq p_1<\dotsm<p_n\leq k}\fT_{-t,a_{p_1}-x_i}\bar\P_{b_{p_1}}e^{(a_{p_2}-a_{p_1})\p^2}\bar\P_{b_{p_2}}\dotsm e^{(a_{p_n}-a_{p_{n-1}})\p^2}\bar\P_{b_{p_n}}\fT_{t,x_j-a_{p_n}}.
\end{equation}
Each summand can be factored as
\begin{equation}
\big(\fT_{-t,a_{p_1}-x_i}\bar\P_{b_{p_1}}\fT_{t,0}\big)\big(\fT_{-t,a_{p_2}-a_{p_1}}\bar\P_{b_{p_2}}\fT_{t,0}\big)\dotsm\big(\fT_{-t,a_{p_n}-a_{p_{n-1}}}\bar\P_{b_{p_n}}\fT_{t,x_j-a_{p_n}}\big).
\end{equation}
By \eqref{eq:Knw0}, as $t\to0$ the first factor goes to ${\bf P}^{\uptext{Hit }\mathfrak{d}_{a_{p_1}}^{b_{p_1}}}_{x_i,a_{p_1}}\uno{x_i\leq a_{p_1}}$, the last factor goes to ${\bf P}^{\uptext{Hit }\mathfrak{d}_{a_{p_n}}^{b_{p_n}}}_{a_{p_{n-1}},x_j}\uno{a_{p_n}\leq x_j}$, and each of the inner factors goes to ${\bf P}^{\uptext{Hit }\mathfrak{d}_{a_{p_s}}^{b_{p_s}}}_{a_{p_{s-1}},a_{p_s}}$, $2\leq s\leq n-1$.
Therefore
\begin{equation}
\lim_{t\to0}e^{-x_i\p^2}\fK^{\hypo(\mathfrak{d}_{\vec a}^{\vec b})}_te^{x_j\p^2}
=\!\!\sum_{n=1}^k(-1)^{n+1}\sum_{1\leq p_1<\dotsm<p_n\leq k}{\bf P}^{\uptext{Hit }\mathfrak{d}_{a_{p_1}}^{b_{p_1}}}_{x_i,a_{p_1}}{\bf P}^{\uptext{Hit }\mathfrak{d}_{a_{p_2}}^{b_{p_2}}}_{a_{p_1},a_{p_2}}\dotsm{\bf P}^{\uptext{Hit }\mathfrak{d}_{a_{p_n}}^{b_{p_n}}}_{a_{p_{n-1}},x_j}\uno{x_i\leq a_{p_1},\,x_j\geq a_{p_n}}
\end{equation}
and then, using inclusion-exclusion again, we deduce finally that
\begin{equation}\label{eq:Kmultinw0}
\lim_{t\to0}e^{-x_i\p^2}\fK^{\hypo(\mathfrak{d}_{\vec a}^{\vec b})}_te^{x_j\p^2}={\bf P}^{\uptext{Hit }\mathfrak{d}_{\vec a}^{\vec b}}_{x_i,x_j}\uno{x_i\leq x_j}
\end{equation}
in operator norm in $L^2([r,\infty))$.

\subsection{Matrix KP initial data}

Now we proceed formally.
Consider compactly supported initial data $\fh\in\UC$, meaning that $\fh(y)=-\infty$ for $y$ outside some compact interval.
Approximating $\fh$ by initial data of the form $\mathfrak{d}_{\vec a}^{\vec b}$ we obtain
\begin{equation}\label{eq:Kmulti}
\lim_{t\to0}e^{-x_i\p^2}\fK^{\hypo(\fh)}_te^{x_j\p^2}={\bf P}^{\uptext{Hit }\fh}_{x_i,x_j}\uno{x_i\leq x_j}.
\end{equation}
In terms of the extended Brownian scattering operator \eqref{eq:scatt-ext}, this gives
\begin{equation}\label{eq:K0}
(\fK_0)_{ij}\coloneqq\fK^{\hypo(\fh)}_{0,\uptext{ext}}(x_i,\cdot;x_j,\cdot)\coloneqq\lim_{t\to0}e^{-x_i\p^2}\fK^{\hypo(\fh  )}_{t,\uptext{ext}}(x_i,\cdot;x_j,\cdot)e^{x_j\p^2}=\begin{dcases*}
-{\bf P}^{\uptext{No hit }\fh}_{x_i,x_j} & if $i<j$,\\
\bar\P_{\fh(x_i)} & if $i=j$,\\
0 & if $i>j$.
\end{dcases*}
\end{equation}

\begin{rem}
This formula formally recovers the correct KPZ fixed point initial data: the kernel $\fK_{ij}(u,v)$ from \eqref{eq:KKscatt} as $t\to0$ becomes $(\fK_0)_{ij}(u+r_i,v+r_j)$ and, in particular, it is upper triangular, so $\det(\fI-\fK)_{L^2(\rr^+)\oplus\dotsm\oplus L^2(\rr^+)}$ as $t\to0$ becomes
\[\det\!\big(\fI-\fK_0\big)_{L^2([r_1,\infty))\oplus\dotsm\oplus L^2([r_n,\infty))}=\prod_{i=1}^n\det\!\big(\fI-\bar\P_{\fh(x_i)}\big)_{L^2([r_i,\infty))}=\prod_{i=1}^n\uno{r_i\geq\fh(x_i)}\]
as desired.
\end{rem}

Now we compute $[\fR_0\fK_0]$ with $\fK_0$ as in \eqref{eq:K0}.
In Sec. \ref{sec:proof} our kernels acted on the $n$-fold direct sum of $L^2(\rr^+)$, but here $\fK_0$ acts on $L^2([r_1,\infty))\oplus\dotsm\oplus L^2([r_n,\infty))$; this corresponds in the current setting to evaluating the $(i,j)$ entry of $(\fI-\fK_0)^{-1}\fK_0$ at $(r_i,r_j)$.
Additionally, it is on this last space that we need to compute compositions of operators; to make this explicit it is convenient to replace the kernel entries $(\fK_0)_{ij}$ by $\P_{r_i}(\fK_0)_{ij}\P_{r_j}$ and simply compute on the $n$-fold direct sum of $L^2(\rr)$.
Doing this, and since $\fK_0$ is upper triangular, we can expand (formally) the entries of $(\fI-\fK_0)^{-1}\fK_0$ as
\begin{equation}
\big((\fI-\fK_0)^{-1}\fK_0\big)_{ij}=\uno{i\leq j}\sum_{\substack{\pi:i\to j\\\pi\uptext{ incr.}}}\prod_{\ell=1}^{|\pi|-1}\P_{r_{\pi(\ell)}}(\fK_0)_{\pi(\ell),\pi(\ell+1)}\P_{r_{\pi(\ell+1)}},\label{eq:RKmultinw-expansion}
\end{equation}
where the sum is over non-decreasing paths $\pi$ going from $i$ to $j$ along integers and $|\pi|$ denotes the length of the path.
Fix $i\leq j$ and assume first that $r_\ell\geq\fh(x_\ell)$ for each $i\leq\ell\leq j$.
Consider a fixed path $\pi$ from $i$ to $j$.
If $\pi(\ell)=\pi(\ell+1)$ for some $\ell$ then the corresponding factor in the product inside the sum will be $\P_{r_{\pi(\ell)}}\bar\P_{\fh(x_{\pi(\ell)})}\P_{r_{\pi(\ell)}}=0$, so only strictly increasing paths contribute to the sum and we get (note that this sum is now finite)
\[\big((\fI-\fK_0)^{-1}\fK_0\big)_{ij}=\uno{i<j}\sum_{\substack{\pi:i\to j\\\pi\uptext{ str. incr.}}}(-1)^{|\pi|-1}\prod_{\ell=1}^{|\pi|-1}\P_{r_{\pi(\ell)}}{\bf P}^{\uptext{No hit } \fh}_{x_{\pi(\ell)},x_{\pi(\ell+1)}}\P_{r_{\pi(\ell+1)}}.\]
Evaluating at $(r_i,r_j)$ 
and applying inclusion-exclusion again, we deduce that as desired (compare with \eqref{eq:matrix0}) that
\begin{equation}\label{eq:abobevo}
\begin{aligned}
[\fR_0\fK_0]_{ij}
&=-\uno{i<j}\pp_{\fB(x_i)=r_i}\!\big(\fB(y)\geq \fh(y)~\forall\,y\in[x_i,x_j],\\
&\hspace{1.6in}\fB(x_\ell)\leq r_\ell\uptext{ for each }x_\ell\in(x_i,x_j),\,\fB(x_j)\in dr_j\big)/dr_j\\
&=-\uno{i<j}{\bf P}^{\geq\fh,\leq-\mathfrak{d}^{-\vec r}_{\vec x}}_{x_i,x_j}(r_i,r_j).
\end{aligned}
\end{equation}

Suppose next that $r_m<\fh(x_m)$ for some $i\leq m\leq j$, and for simplicity assume that this is the only such index satisfying the condition (the argument can be generalized easily).
Assume also that $i<j$.
From the argument in the previous case we know that if $\pi\!:i\to j$ has a constant piece which stays at any index other than $m$, then $\pi$ does not contribute to the sum in \eqref{eq:RKmultinw-expansion}.
Hence any path $\pi$ from $i$ to $j$ which does contribute to the sum can be decomposed as $\pi_1\circ\upsilon\circ\pi_2$ with $\pi_1\!:i\to m$ and $\pi_2\!:m\to j$ strictly increasing (we allow for $\pi_1$ or $\pi_2$ to be empty if $m=i$ or $m=j$), and $\upsilon$ staying at $m$ for a given number of steps (which could be $0$).
The product inside the sum in \eqref{eq:RKmultinw-expansion} splits between factors coming from the three pieces of the path, and from the middle part we get a factor $\big(\P_{r_\ell}\bar\P_{\fh(x_m)}\P_{r_m}\big)^{|\upsilon|}=\fI\cdot\uno{|\upsilon|=0}+\P_{r_m}\bar\P_{\fh(x_m)}\uno{|\upsilon|>0}$.
In other words, and repeating the previous argument,
\begin{align}
&\big((\fI-\fK_0)^{-1}\fK_0\big)_{ij}=-\sum_{\substack{\pi_1:i\to m\\\pi_1\uptext{ str. incr.}}}\sum_{\substack{\pi_2:m\to j\\\pi_2\uptext{ str. incr.}}}\sum_{\nu\geq0}\left((-1)^{|\pi_1|-1}\prod_{\ell=1}^{|\pi_1|-1}\P_{r_{\pi_1(\ell)}}{\bf P}^{\uptext{No hit } \fh}_{x_{\pi_1(\ell)},x_{\pi_1(\ell+1)}}\P_{r_{\pi_1(\ell+1)}}\right)\\
&\hspace{0.75in}\times\big(\fI\cdot\uno{\nu=0}+\P_{r_m}\bar\P_{\fh(x_m)}\uno{\nu>0}\big)\left((-1)^{|\pi_2|-1}\prod_{\ell=1}^{|\pi_2|-1}\P_{r_{\pi_2(\ell)}}{\bf P}^{\uptext{No hit } \fh}_{x_{\pi_2(\ell)},x_{\pi_2(\ell+1)}}\P_{r_{\pi_2(\ell+1)}}\right)\\
&\qquad=-\sum_{\nu\geq0}\P_{r_i}{\bf P}^{{\geq\fh,\leq\mathfrak{d}_{\vec x}^{\vec r}}}_{x_i,x_m}\big(\fI\cdot\uno{\nu=0}+\P_{r_m}\bar\P_{\fh(x_m)}\uno{\nu>0}\big){\bf P}^{{\geq\fh,\leq\mathfrak{d}_{\vec x}^{\vec r}}}_{x_m,x_j}\P_{r_j}.
\end{align}
But ${\bf P}^{{\geq\fh,\leq\mathfrak{d}_{\vec x}^{\vec r}}}_{x_i,x_m}(u,v)=0$, because at the endpoint $v$ it requires $\fh(x_m)\leq v\leq r_m$ (the analogous statement holds for the other factor).
Hence we conclude that, in this case, $[\fR_0\fK_0]_{i,j}=0$, which for the same reason means that \eqref{eq:abobevo} still holds.

Suppose finally that $i=j$ and $r_i<\fh(x_i)$.
Now the only possible paths in \eqref{eq:RKmultinw-expansion} are constant paths $\pi$ of arbitrary length $|\pi|\geq1$.
Each such $|\pi|$ contributes a term of the form $\P_{r_i}\bar\P_{\fh(x_i)}$, which evaluated at $(r_i,r_i)$ is taken to be $1$, and hence $[\fR_0\fK_0]_{i,i}$ diverges to $\infty$ in this case (which coincides with the physical meaning of this quantity, namely $\p_{r_i}\tsm\log F(t,x_1,\dotsc,x_n,r_1,\dotsc,r_n)$).

The conclusion is then that $Q(0,x_1,\dotsc,x_n,r_1,\dotsc,r_n)=[\fR_0\fK_0]$ satisfies \eqref{eq:matKPini}.

\section{Alternative derivation of KP-II for narrow wedge multipoint distributions}
\label{pif}

In this section we will derive the KP-II equation \eqref{eq:KP-II-Airy} for the Airy$_2$ process directly using the path-integral formula for the KPZ fixed point \cite[Prop. 4.3]{fixedpt}.
Define $F(t,\vec x+y,\vec r+a)$ as in the first equality of \eqref{eq:sktime-nw}.
Then letting $\fK_{t,x}=\fK^{\hypo(\fh_0)}_{t}(x,\cdot;x,\cdot)$ we have
\begin{align}
F(t,\vec x+y,\vec r+a)&=\det(\fI-\fK_{t,x_1+y}+\bar\P_{r_1+a}e^{(x_2-x_1)\p^2}\bar\P_{r_2+a}\dotsm\bar\P_{r_m+a}e^{(x_1-x_m)\p^2}\fK_{t,x_1+y})\\
&=\det(\fI-\fK+\bar\P_{r_1}e^{(x_2-x_1)\p^2}\bar\P_{r_2}\dotsm\bar\P_{r_m}e^{(x_1-x_m)\p^2}\fK)
\end{align}
with $\fK=\fK_{t,x_1+y}(a+\cdot,a+\cdot)=e^{a\p}\fK_{t,x_1+y}e^{-a\p}$.
Note that the product of operators preceding $\fK$ in the last term does not depend on $t$, $y$ or $a$; call it $\fI-\fP$ so that $F=\det(\fI-\fP\fK)$.
Up to here this is general, but now we specialize to the narrow wedge case, for which $\fK=e^{a\p}(\fT_{t,-x_1-y})^*\bar\P_0\fT_{t,x_1+y}e^{-a\p}$ (see \eqref{eq:nwKt}).
Using the cyclic property of the determinant we get
\begin{align}
F&=\det(\fI-\bar\P_0\fT_{t,x_1+y}e^{-a\p}\fP e^{a\p}(\fT_{t,-x_1-y})^*\bar\P_0)=\det(\fI-\bar\P_0e^{-a\p}\fT_{t,x_1+y}\fP (\fT_{t,-x_1-y})^*e^{a\p}\bar\P_0)\\
&=\det(\fI-\P_0\varrho e^{-a\p}\fT_{t,x_1+y}\fP (\fT_{t,-x_1-y})^*e^{a\p}\varrho\P_0)
\end{align}
with $\varrho$ the reflection operator $\varrho f(x)=f(-x)$. 
So letting
\[\fL=\varrho e^{-a\p}\fT_{t,x_1+y}\fP (\fT_{t,-x_1-y})^*e^{a\p}\varrho=e^{a\p}(\fT_{t,x_1+y})^*(\varrho\fP\varrho)\fT_{t,-x_1-y}e^{-a\p}\]
(the second equality is a simple computation) we get that $F=\det(\fI-\fL)$.
Now $\p_a\fL(u,v)=(\p_u+\p_v)\fL(u,v)$, $\p_t\fL(u,v)=-\frac13(\p_u^3+\p_v^3)\fL(u,v)$, and $\p_y\fL(u,v)=(\p_u^2-\p_v^2)\fL(u,v)$, which correspond to (1), (2) and (3) of Thm. \ref{thm3} (except for the change $y\mapsto-y$, which as in Sec. \ref{sec:airy} is irrelevant), so the theorem implies that $\phi=\p_a^2\log(F)$ solves KP-II in $(t,y,a)$, and translating back to the $\cD_r$, $\cD_x$ derivatives yields \eqref{eq:KP-II-Airy}.

\vs

\noindent{\bf Acknowledgements.}
The authors thank P. Le Doussal for pointing out that a claim made in the first version of the present paper about the generating function of the KPZ equation with flat initial condition solving KdV was inconsistent with known moment formulas \cite{cal-led}.
JQ was supported by the Natural Sciences and Engineering Research Council of Canada and by a Killam research fellowship.
DR was supported by CMM ANID PIA AFB170001, by Programa Iniciativa Cient\'ifica Milenio grant number NC120062 through Nucleus Millenium Stochastic Models of Complex and Disordered Systems, and by Fondecyt Grant 1201914.

\printbibliography[heading=apa]

\end{document}